\documentclass[10pt, twoside, openany]{book}

  \usepackage{amsmath, amsthm, amssymb, slashed, framed}
\usepackage{titlesec}

 \usepackage {graphicx}
\usepackage{amsfonts}
\usepackage{amsthm}
\usepackage{amsmath}
\usepackage{latexsym}
\usepackage{amssymb}
\usepackage{epsfig,color}

\usepackage{graphicx, amssymb}
 \usepackage{mathtools}
  
\newtheorem{definition}{Definition}[chapter]
\newtheorem{exercise}[]{Exercise}

\newtheorem{example}[definition]{Example} 
\newtheorem{theorem}[definition]{Theorem} 
\newtheorem{remark}[definition]{Remark} 
\newtheorem{proposition}[definition]{Proposition}

\newtheorem{rem}[definition]{Remark}

\newcommand{\Na}{\begin{eqnarray}}
\newcommand{\ENa}{\end{eqnarray}}

\usepackage[english]{}

\usepackage[a4paper,inner=30mm,outer=30mm]{geometry}

\titleformat{\chapter}[display]
  {\bfseries\Large}
  {\filright\MakeUppercase{\chaptertitlename} \Huge\thechapter}
  {1ex}
  {\titlerule\vspace{1ex}\filleft}
  [\vspace{1ex}\titlerule]

\makeatletter
\def\maketitle{%
  \null
  \thispagestyle{empty}%
  \vfill
  \begin{center}\leavevmode
    \normalfont
    {\LARGE\raggedleft \@author\par}%
    \hrulefill\par
    {\huge\raggedright \@title\par}%
    \vskip 1cm
  \end{center}%
  \vfill
  \null
  \cleardoublepage
  }
\makeatother
\author{Laura P. Schaposnik}
\author{Laura P. Schaposnik}
\title{ An introduction to spectral data \\for Higgs bundles}
\date{}

\begin{document}

\let\cleardoublepage\clearpage

\maketitle

\frontmatter

%
%
%
%
%
%
%
%
%
%
%
%
 
\newpage

\let\cleardoublepage\clearpage
\mainmatter
\sloppy


\chapter*{Outline}
\begin{small}
These notes have been prepared as reading material for the mini-course ``{\em An introduction to spectral data for Higgs bundles}'' that the author gave at the National University of Singapore as part of the  {\em Summer school on the moduli space of Higgs bundles}.  \end{small}
\subsection*{Lecture \ref{cap1}}
\begin{small}
The first lecture shall introduce classical Higgs bundles and the Hitchin fibration, and describe the associated spectral data in the case of principal Higgs bundles for classical complex Lie groups. Whilst  bibliography is provided in the text, the main references followed are Hitchin's papers \cite{N1, N2, N5, N3}.
\end{small}
\subsection*{Lecture \ref{cap2}}
\begin{small}
During the second lecture we shall construct Higgs bundles for real forms of classical complex Lie groups as fixed points of involutions, and describe the corresponding spectral data when known, as appearing in  \cite{Lau, thesis, umm} and \cite{non ab}. Along the way, we shall mention different applications and open problems related to the methods introduced in both lectures. 
\end{small}
\subsubsection*{Exercises}

\begin{small}
 Each lecture contains exercises of varying difficulty, whose solutions can be found in \cite{thesis}, and commented in the .tex file. Open problems  which might be tackled with methods similar to the ones introduced in the lectures shall also be mentioned, and appear indicated with \textbf{((*))}.  For each of these problems we  suggest references which feature approaches that may be useful.  \end{small}
\subsubsection*{Bibliography}
\begin{small}
\noindent We shall highlight the main references considered, as well as the precise places where the methods used were developed. Since it proves to be very difficult to give a comprehensive and exhaustive account of research in tangential areas, we shall restrain ourselves to mentioning related work only when it directly involves methods using {\em spectral data}. The reader should  refer to references in the bibliography for further research in related topics (e.g., see references in \cite{ap, ana1, thesis}). \end{small}

\subsubsection*{Acknowledgments}
\begin{small}
\noindent   The author would  like to thank the Institute of Mathematical Sciences at the National University of Singapore and the organisers of the program {\it The Geometry, Topology and Physics of Moduli Spaces of Higgs Bundles} for providing an ideal environment for research and collaborations. The author also acknowledges support from U.S. National Science Foundation grants DMS 1107452, 1107263, 1107367 ``sRNMS: GEometric structures And Representation varieties" (the GEAR Network). \end{small}

\chapter{Spectral data for $G_c$-Higgs bundles}\label{cap1}

 \hfill{\vbox{\hbox{\it  The art of doing mathematics } \hbox{\it  consists in finding that special} \hbox{\it case which contains all the  germs }\hbox{\it of generality. } \hbox{\it\noindent\rule{6cm}{0.4pt}}}}
 
 \hfill{\vbox{\hbox{{David Hilbert}}}

\bigskip

Following \cite{N1,N2, N3} we dedicate this lecture to overview classical Higgs bundles  as well as $G_{c}$-Higgs bundles for the  groups $G_{c}= SL(n,\mathbb{C})$, $Sp(n,\mathbb{C})$, $SO(2n+1,\mathbb{C})$ and $SO(2n,\mathbb{C})$. In each case we introduce the Hitchin fibration and describe the generic fibres  through \textit{spectral data}, i.e., an associated spectral curve and a line bundle on it.

\section{$G_c$-Higgs bundles}  \label{sec:mod}

 Consider $\Sigma$ a compact Riemann surface of genus $g\geq 2$ with canonical bundle $K=T^*\Sigma$. Classically, a Higgs bundle  on $\Sigma$ is defined as follows:

  \begin{definition}\label{def:clasical}
   A {\em Higgs bundle} is a pair $(E,\Phi)$ for $E$ a holomorphic vector bundle on $\Sigma$, and $\Phi$, the {\em  Higgs field}, a holomorphic section in $H^{0}(\Sigma,\text{End}(E)\otimes K)$.    \end{definition}
In order to understand better what Higgs bundles are and how to generalise the definition, we shall first look at the moduli space of vector bundles and then study the moduli space of classical Higgs bundles and its associated spectral data. For more details  the reader should refer to \cite{N1,N2,donald,cor,simpson88,nit,simpson}.

\subsection{Moduli space of vector bundles}

Holomorphic vector bundles $E$ on a compact Riemann surface $\Sigma$ of genus $g\geq 2$ are  topologically classified by their rank $rk(E)$  and their  degree $deg(E)$.
\begin{definition}
 The {\em  slope} of a holomorphic vector bundle $E$  is  
$\mu(E):= deg(E)/  rk(E)$ and is used to define stability conditions:
  A vector bundle $E$ is said to be {\em stable} ({\em semi-stable}) if for any proper, non-zero sub-bundle $F\subset E$ one has $\mu(F)<\mu(E)$ ($\mu(F)\leq \mu(E)$).
 It is {\em  polystable} if it is a direct sum of stable bundles whose slope is the same  as $E$.
  \end{definition}

 It is known that the space of holomorphic bundles of fixed rank and fixed degree, up to isomorphism, is not a Hausdorff space. However, through Mumford's Geometric Invariant Theory one can construct the moduli space $\mathcal{N}(n,d)$ of stable bundles of fixed rank $n$ and degree $d$, which has the natural structure of an algebraic variety.

\begin{remark}
For coprime $n$ and $d$, the moduli space $\mathcal{N}(n,d)$  is a smooth projective algebraic variety of dimension $n^{2}(g-1)+1$.
\end{remark}

\begin{remark}
All line bundles are stable, and thus $\mathcal{N}(1,d)$ contains all line bundles of degree $d$, and is isomorphic to $\text{Jac}^{d}(\Sigma)$ of $\Sigma$, an abelian variety of dimension $g$.
\end{remark}

Let $G_{c}$ be a complex semisimple Lie group.  Following \cite{ram} one can define stability for principal $G_{c}$-bundles as follows (see \cite[Section 1.1]{ap} for a comprehensive study):

\begin{definition}
 A holomorphic principal $G_{c}$-bundle $P$ is said to be {\em stable} ({\em semi-stable}) if  for every reduction $\sigma : \Sigma \rightarrow P/Q $ to maximal parabolic subgroups $Q$ of $G_{c}$ one has $\text{deg}~ \sigma^{*} T_{rel}>0 ~ ( \geq 0),$
where $T_{rel}$ is the relative tangent bundle for the projection $P/Q\rightarrow \Sigma$.
\end{definition}

The notion of polystability may be carried over to principal $G_{c}$-bundles, allowing one to construct the moduli space of isomorphism classes of polystable principal $G_{c}$-bundles of fixed topological type over the compact Riemann surface $\Sigma$.

\subsection{Moduli space of classical Higgs bundles}

In order to define the moduli space of Higgs bundles, the following stability condition is considered:
\begin{definition}
 A vector subbundle $F$ of $E$ for which $\Phi(F)\subset F\otimes K$ is said to be a $\Phi$-{\em invariant subbundle} of $E$.  
 A Higgs bundle $(E,\Phi)$ is said to be
\begin{itemize}
 \item {\em  stable} ({\em semi-stable}) if for each proper $\Phi$-invariant  $F\subset E$ one has $\mu(F)<~\mu(E)~(equiv. \leq)$;
  \item {\em   polystable} if $(E,\Phi)=(E_{1},\Phi_{1})\oplus (E_{2},\Phi_{2})\oplus \ldots \oplus (E_{r},\Phi_{r})$, where $(E_{i},\Phi_{i})$ is stable with $\mu (E_{i})=\mu(E)$ for all $i$.
\end{itemize}

\end{definition}
\begin{remark}\label{invariantsubbundle}
 The characteristic polynomial of $\Phi$ restricted to an invariant subbundle  divides the characteristic polynomial of $\Phi$. 
\end{remark}

\begin{example} \label{exa}\label{classic example} Choose a square root $K^{1/2}$ of the canonical bundle $K$, and  a section $\omega$ of $K^{2}$. A family of classical Higgs bundles $(E,\Phi_\omega)$ may be obtained by considering the vector bundle  $E=K^{\frac{1}{2}}\oplus K^{-\frac{1}{2}}$ and the Higgs bundle   $\Phi_{\omega}$  given by
\begin{eqnarray}
 \Phi_{\omega}=\left(\begin{array}{cc}
0&\omega\\1&0
            \end{array}\right) \in H^{0}(\Sigma, \text{End}(E)\otimes K)
\nonumber
\end{eqnarray} 
\end{example}
\begin{exercise} Show that the pairs $(E, \Phi_\omega)$ from Example \ref{classic example} are stable. 
\end{exercise}

 

\begin{exercise}Prove that if a Higgs bundle $(E,\Phi)$ is stable, then for $\lambda\in \mathbb{C}^{*}$ and $\alpha$ a holomorphic automorphism of $E$, the induced Higgs bundles $(E,\lambda \Phi)$ and $(E,\alpha^{*}\Phi)$ are stable.
\end{exercise}


In order to define the moduli space of classical Higgs bundles, we shall first define an appropriate equivalence relation. For this, consider a strictly semi-stable Higgs bundle $(E, \Phi)$. As it is not stable, $E$ admits a subbundle $F\subset E$ of the same slope which is preserved by $\Phi$. If $F$ is a subbundle of $E$ of least rank and same slope which is preserved by $\Phi$, it follows that $F$ is stable and hence the induced pair $(F,\Phi)$ is stable.  Then, by induction one obtains a flag of subbundles
$F_{0}=0\subset F_{1}\subset \ldots \subset F_{r}=E$
where $\mu(F_{i}/F_{i-1})=\mu(E)$ for $1\leq i\leq r$, and where the induced Higgs bundles $(F_{i}/F_{i-1}, \Phi_{i})$ are stable. This is the \textit{Jordan-H\"{o}lder filtration} of $E$, and it is not unique.  However, the graded object
$\text{Gr}(E,\Phi):=\bigoplus_{i=1}^{r}(F_{i}/F_{i=1},\Phi_{i})$
is unique up to isomorphism.

\begin{definition}
 Two semi-stable Higgs bundles $(E,\Phi)$ and $(E',\Phi')$ are said to be $S${\em -equivalent } if $\text{Gr}(E,\Phi)\cong \text{Gr}(E',\Phi')$. 
\end{definition}

\begin{exercise} If a pair $(E,\Phi)$ is strictly stable, what is the induced Jordan-H\"{o}lder filtration? 
\end{exercise}

Following   \cite[Theorem 5.10]{nit} we let $\mathcal{M}(n,d)$ be the moduli space of $S$-equivalence classes of semi-stable Higgs bundles of fixed degree $d$ and fixed rank $n$. The moduli space $\mathcal{M}(n,d)$  is a quasi-projective scheme, and has an open subscheme $\mathcal{M}'(n,d)$ which is the moduli scheme of stable pairs. Thus, every point is represented by either a stable or a polystable Higgs bundle. When $d$ and $n$ are coprime, the moduli space $\mathcal{M}(n,d)$  is smooth.

The cotangent space of $\mathcal{N}(n,d)$ over the stable locus is contained in $\mathcal{M}(n,d)$ as a Zariski  open subset.  The moduli space $\mathcal{M}(n,d)$  is a non-compact  variety which has complex dimension $ 2n^{2}(g-1)+2$. Moreover,  it is a hyperk\"ahler manifold with natural symplectic form $\omega$ defined on the infinitesimal deformations
 $(\dot A,\dot \Phi)$ of a Higgs bundle $(E,\Phi)$   by
\begin{equation}\label{ch2:2.1}
 \omega((\dot{A}_{1},\dot{\Phi}_{1}),(\dot{A}_{2},\dot{\Phi}_{2}))=\int_{\Sigma}\text{tr}(\dot{A}_{1}\dot{\Phi}_{2}-\dot{A}_{2}\dot{\Phi}_{1}),
\end{equation}
 where $\dot A \in \Omega^{0,1}(\text{End}_{0} E)$ and $\dot\Phi\in \Omega^{1,0}(\text{End}_0 E)$ 
(see \cite{N1, N2} for details). For simplicity, we shall fix $n$ and $d$ and write $\mathcal{M}$ for $\mathcal{M}(n,d)$.

\subsection{Moduli space of $G_c$-Higgs bundles}
 
The notion of Higgs bundle can be generalized to encompass principal $G_c$-bundles, for $G_c$ a complex  semi-simple Lie group. For more details, the reader should refer  to \cite{N2}.
  \begin{definition}\label{principalLie}\label{defHiggs}
   A $G_{c}${\em -Higgs bundle} is a pair $(P,\Phi)$ where  $P$ is a principal $G_{c}$-bundle over $\Sigma$, and the Higgs field $\Phi$ is a holomorphic section of the vector bundle  {\em ad}$P\otimes_{\mathbb{C}} K$, for {\em ad}$P$ the vector bundle associated to the adjoint representation. 
  \end{definition}

When $G_{c}\subset GL(n,\mathbb{C})$, a $G_{c}$-Higgs bundle gives rise to a Higgs bundle in the classical sense,  with
some extra structure reflecting the definition of $G_{c}$. In particular, classical Higgs bundles are given by $GL(n,\mathbb{C})$-Higgs bundles.  
\begin{example} The Higgs bundles in Example \ref{exa} have traceless Higgs field, and the determinant  $\Lambda^{2}E$ is trivial. Hence, for each quadratic differential $\omega$ one has an $SL(2,\mathbb{C})$-Higgs bundle $(E,\Phi_\omega)$.
\end{example}
By extending the stability definitions for principal $G_{c}$-bundles, one can define \textit{stable}, \textit{semi-stable} and $polystable$ $G_{c}$-Higgs bundles. Since in these notes we shall be working with Higgs pairs which do not preserve any subbundle, they will be automatically stable and thus we shall not dedicate time to recall the main study of stability for principal Higgs bundles. For details about the corresponding constructions, the reader should refer for example to \cite[Section 3]{biswas}, \cite[Section 1]{ap}. 
We denote by $\mathcal{M}_{G_{c}}$  the moduli space of $S$-equivalence classes of   polystable $G_{c}$-Higgs bundles.

In the remainder of this Section, following \cite{N2} and \cite{N3} we  introduce the Hitchin fibration and describe the generic fibres  for  $G_c$-Higgs bundles where $G_{c} = GL(n,\mathbb{C})$, $SL(n,\mathbb{C})$, $Sp(2n,\mathbb{C})$, $SO(2n+1,\mathbb{C})$ and $SO(2n,\mathbb{C})$. We shall cover with more detail the initial cases, and leave as an exercise to the reader some of the results for the latter groups. 

\section{The Hitchin fibration}
A natural way of studying $\mathcal{M}_{G_c}$ is through the Hitchin fibration, as introduced in \cite{N2}. 
We shall denote by $p_{i}$,  for $i=1,\ldots, k$, a homogeneous basis for the algebra of invariant polynomials of the Lie algebra  $\mathfrak{g}_{c}$ of $G_{c}$, and let $d_{i}$ be their degrees. Then, the {\em Hitchin fibration} is given by
\begin{eqnarray}
 h~:~ \mathcal{M}_{G_{c}}&\longrightarrow&\mathcal{A}_{G_{c}}:=\bigoplus_{i=1}^{k}H^{0}(\Sigma,K^{d_{i}}),\\
 (E,\Phi)&\mapsto& (p_{1}(\Phi), \ldots, p_{k}(\Phi)).
\end{eqnarray}
The map $h$ is referred to as the {\em Hitchin~map}, and is a proper map for any choice of basis (see \cite[Section 4]{N2} for details). Furthermore,  $\mathcal{A}_{G_{c}}$ always satisfies
$\text{dim} \mathcal{A}_{G_{c}} =\text{dim}\mathcal{M}_{G_{c}}/2,$
making the Higgs bundle moduli space into an integrable system.
%


\begin{remark}
A homogenous basis of invariant polynomials for classical Higgs bundles $(E,\Phi)$ of rank $n$ can be taken as {\em tr}$(\Phi^i)$ for $1\leq i\leq n$.
\end{remark}

\begin{remark} Whilst a formal   definition of the smooth locus of the Hitchin base can be given (e.g., see \cite{donagi}) in these lectures we shall note that the  generic fibres of the Hitchin fibration are smooth, and thus generic points in the Hitchin base are in the smooth locus.  
 \end{remark}

\subsection{$GL(n,\mathbb{C})$-Higgs bundles}

As before, let $K$ be the canonical bundle of $\Sigma$, and $X$ its total space with projection $\pi:X\rightarrow \Sigma$. We shall denote by $\eta$ the tautological section of the pull back $\pi^{*}K$ on $X$. Abusing notation we denote with the same symbols the sections  of powers $K^{i}$ on $\Sigma$ and their pull backs to $X$. The characteristic polynomial of a Higgs bundle $(E,\Phi)$ in a generic fibre $h^{-1}(a)$ defines a smooth curve $\pi:S_a\rightarrow \Sigma$ in $X$, the {\em spectral curve} of $\Phi$, whose equation is
\begin{eqnarray}
\text{det}(\eta Id- \pi^*\Phi)= \eta^{n}+a_{1}\eta^{n-1}+a_{2}\eta^{n-2}+\ldots + a_{n-1}\eta+a_{n}=0,\label{smootheq}
\end{eqnarray}
for $a_{i}\in H^{0}(\Sigma,K^{i})$ (for simplicity, we shall write $\text{det}(\eta-\Phi)$ for the characteristic polynomial of the Higgs field $\Phi$, and drop the subscript $a$ of $S_a$).  By the adjunction formula on $X$ (see e.g. \cite{harris}), since the canonical bundle $K$ has trivial cotangent bundle one has $K_{S}\cong \pi^{*}K^{n}$, and hence the genus of $S$ is
 \begin{eqnarray}\label{classic curve genus}g_{S}=1+n^{2}(g-1).\label{genusSclassical}\end{eqnarray}
%
 
 \begin{framed} The {\em spectral data} classical Higgs bundles in a smooth fibre   of the Hitchin fibration  is given by a spectral curve $S$ defined as in (\ref{smootheq}) and a line bundle $L \in$ Jac$(S)$. \end{framed}
 
 In order to see that the smooth fibres of the Hitchin fibration are Jacobians, 
 starting with a line bundle $L$ on the smooth curve $\pi:S\rightarrow \Sigma$ with equation as in (\ref{smootheq}), we shall obtain a classical Higgs bundle   by considering the direct image $\pi_{*}L$ of $L$. Recall that by definition of direct image, given an open set $\mathcal{U}\subset \Sigma$, one has
$  H^{0}(\pi^{-1}(\mathcal{U}),L)= H^{0}(\mathcal{U}, \pi_{*}L).$
 Multiplication by the the tautological section $\eta$ induces the map
\begin{eqnarray}
 H^{0}(\pi^{-1}(\mathcal{U}),L)\xrightarrow{\eta} H^{0}(\pi^{-1}(\mathcal{U}),L\otimes\pi^{*}K),\nonumber
\end{eqnarray}
which by definition of direct image can be pushed down to give \begin{eqnarray}
 \Phi:    \pi_{*}L\rightarrow  \pi_{*}L \otimes K.\nonumber
\end{eqnarray}
Then, one obtains a Higgs field $\Phi\in H^{0}(\Sigma,\text{End}E \otimes K)$ for $E:=\pi_{*}L$.

\begin{exercise}
 Use Grothendieck-Riemann-Roch to see that $\text{deg}(E)=\text{deg}(L) + (n^{2}-n)(1 - g).$ 
 \end{exercise}Moreover, the Higgs field satisfies its characteristic equation, which by construction is given by 
$\eta^{n}+a_{1}\eta^{n-1}+a_{2}\eta^{n-2}+\ldots + a_{n-1}\eta+a_{n}=0.$
Furthermore, since $S$ is irreducible, from Remark \ref{invariantsubbundle} there are no invariant subbundles of the Higgs field, making the induced Higgs bundle $(E,\Phi)$  stable.

%
%


Conversely, let $(E,\Phi)$ be a classical Higgs bundle. The characteristic polynomial is given by
$\text{det}(x -\Phi)=x^{n}+a_{1}x^{n-1}+a_{2}x^{n-2}+\ldots + a_{n-1}x+a_{n},$
and its coefficients  define the $spectral~curve$ $S$ in the total space $X$ whose equation is (\ref{smootheq}).

From \cite[Proposition 3.6]{bobi}, there is a bijective correspondence between Higgs bundles $(E,\Phi)$ and the line bundles $L$ on the spectral curve $S$ described previously. This correspondence identifies the fibre of the Hitchin map with the Picard variety of line bundles of the appropriate degree.  By tensoring the line bundles $L$ with a chosen line bundle of degree $-\text{deg}(L)$, one obtains a point in the Jacobian $\text{Jac}(S)$, the abelian variety of line bundles of  degree zero on $S$, which has dimension $g_{S}$ as in (\ref{genusSclassical}). In particular, the Jacobian variety is the connected component of the identity in the Picard group $H^{1}(S,\mathcal{O}^{*}_{S})$.
Thus, the fibre of the classical Hitchin fibration $h :\mathcal{M}\rightarrow \mathcal{A} $ is isomorphic to the Jacobian of the spectral curve $S$. For more details, the reader should refer for example to \cite[Section 2]{N3}.

\begin{example}\label{newexample}
In the case of a classical rank 2 Higgs bundle $(E,\Phi)$, the characteristic polynomial of $\Phi$ defines a spectral curve $\pi:S \rightarrow \Sigma$. This is a 2-fold cover of $\Sigma$ in the total space of $K$, and has  equation
$\eta^{2}+a_{2}=0$,  
for $a_{2}$ a quadratic differential  and $\eta$ the tautological section of $\pi^{*}K$. By \cite[Remark 3.5]{bobi} the curve is smooth when $a_{2}$ has simple zeros, and in this case the ramification points are  given by the divisor of $a_{2}$. For $z$ a  local coordinate near a  ramification point, the covering is given by
$ z\mapsto z^{2}:=w.$ In a neighbourhood of $z=0$, a section of the line bundle $M$ looks like
$f(w)=f_{0}(w)+zf_{1}(w).$
Since the Higgs field is obtained via multiplication by $\eta$, one has
\begin{eqnarray}
\Phi (f_{0}(w)+zf_{1}(w)) = w f_{1}(w)+ z f_{0}(w),
\end{eqnarray}
and thus  a local form of the Higgs field $\Phi$ is given by 
\[\Phi=\left(\begin{array}
              {cc}
0&w\\
1&0
             \end{array}
\right).\]
                              

\end{example}

\begin{remark}
When $G_c\subset GL(n,\mathbb{C})$, the spectral data of a $G_c$-Higgs bundle is given by the spectral data of the pair as a classical Higgs bundle, satisfying extra conditions. 
\end{remark}

\begin{remark}  For generic $G_{c}$, a  description of the fibres can be obtained by means of Cameral covers  \cite{donagi1} (see also \cite{donagi}) which is equivalent to the one given in the next sections for classical Lie groups. 
\end{remark}

\subsection{$SL(n,\mathbb{C})$-Higgs bundles }

When $G_{c}=SL(n,\mathbb{C})$ we apply Definition \ref{defHiggs} to obtain the following:

\begin{definition}An $SL(n,\mathbb{C})$-Higgs bundle is a classical Higgs bundle $(E,\Phi)$ where the rank $n$ vector bundle $E$ has trivial determinant and the Higgs field has zero trace.  
\end{definition}

A basis for the invariant polynomials on the Lie algebra $\mathfrak{sl}(n,\mathbb{C})$ is given by the coefficients of the characteristic polynomial of a trace-free matrix $A\in \mathfrak{sl}(n,\mathbb{C})$.
 In this case, the spectral curve $\pi: S\rightarrow \Sigma$ associated to the Higgs bundle has equation
\begin{eqnarray}\label{spectralsl}
 \eta^{n}+a_{2}\eta^{n-2}+\ldots + a_{n-1}\eta+a_{n}=0,
\end{eqnarray}
where $a_{i}\in H^{0}(\Sigma,K^{i})$ are the coefficients of the characteristic polynomial of the Higgs field $\Phi$. Generically   $S$ is a smooth curve of genus $g_{S}=1+n^{2}(g-1),$ and the coefficients define the corresponding 
 Hitchin fibration
\begin{eqnarray}
 h~:~ \mathcal{M}_{SL(n,\mathbb{C})}\longrightarrow\mathcal{A}_{SL(n,\mathbb{C})}:=\bigoplus_{i=2}^{n}H^{0}(\Sigma,K^{i}).
\end{eqnarray}

In this case the generic fibres of the Hitchin fibration are given by the subset of $\text{Jac}(S)$  of line bundles $L$ on $S$ for which $\pi_{*}L=E$ and $\Lambda^{n}\pi_{*}L$ is trivial. By understanding these conditions in terms of $L$ one has the following:

   \begin{framed}
 The generic fibre $h^{-1}(a)$ of the $SL(n,\mathbb{C})$ Hitchin fibration is biholomorphically equivalent to the Prym$(S,\Sigma)$ variety of the  spectral curve $S_a$ defined as in (\ref{smootheq}).
\end{framed}

In order to see why one has to take the Prym variety, recall that the Norm map
$$\text{Nm}: \text{Pic}(S)\rightarrow \text{Pic}(\Sigma),$$
 associated to $\pi$  is defined on  divisor classes by $\text{Nm}( \sum n_{i}p_{i})=\sum n_{i} \pi(p_{i})$. 
In particular,
\[\text{Nm}(\pi^{-1}(x))= \pi(\pi^{-1}(x)) = nx.\]
The kernel of the Norm map is the $Prym ~variety$, and is denoted by $\text{Prym}(S,\Sigma)$.
 From \cite[Section 4]{bobi}, the determinant bundle of $L$ satisfies
\[\Lambda^{n}\pi_{*}L\cong \text{Nm}(L)\otimes K^{-n(n-1)/2}.\]
 Thus, $\Lambda^{n}\pi_{*}L$ is trivial if and only if
\begin{eqnarray}
  \text{Nm}(L) \cong K^{n(n-1)/2}. \label{ab1}
\end{eqnarray}
Equivalently, since $\text{Nm}( \sum n_{i}\pi^{-1}(p_{i}))= n \sum n_{i}p_{i},$ the determinant bundle $\Lambda^{n}\pi_{*}L$ is trivial
 if the line bundle $M:=L\otimes \pi^{*}K^{-(n-1)/2}$ is in the Prym variety. 
 
 \begin{remark}In the case of even rank, equation (\ref{ab1}) implies a choice of a square root of $K$  (see \cite{N1} and \cite[Section 2.2]{N3} for more details).   
\end{remark}

\subsection{$Sp(2n,\mathbb{C})$-Higgs bundles  } \label{subsec:sp}

Let $G_{c}=Sp(2n,\mathbb{C})$, and let $V$ be a $2n$ dimensional vector space with a non-degenerate
skew-symmetric form $< , >$. For $v_{i},v_{j}$ eigenvectors of $A\in \mathfrak{sp}(2n,\mathbb{C})$ with eigenvalues $\lambda_{i}$ and $\lambda_{j}$,  
\begin{eqnarray}
 \lambda_{i}<v_{i},v_{j}>&=&<\lambda_{i}v_{i},v_{j}> \\
&=&<Av_{i},v_{j}>\\&=&-<v_{i},Av_{j}> \\&=&-<v_{i},\lambda_{j}v_{j}> 
=-\lambda_{j}<v_{i},v_{j}>.\nonumber
\end{eqnarray}
From the above one has that $<v_{i}, v_{j}>=0$ unless $\lambda_{i}=-\lambda_{j}$. Since $<v_{j},v_{j}>=0$, from the non-degeneracy of the symplectic inner product it follows that if $\lambda_{i}$ is an eigenvalue so is $-\lambda_{i}$. Thus, distinct eigenvalues  of $A$ must occur in $\pm \lambda_{i}$ pairs, and the corresponding eigenspaces are paired by the symplectic form. The  characteristic polynomial of $A$ must therefore be of the form
\[\text{det}(x-A)=x^{2n}+a_{1}x^{2n-2}+\ldots+a_{n-1}x^{2}+a_{n},\]
and a basis for the invariant polynomials on the Lie algebra $\mathfrak{sp}(2n,\mathbb{C})$ is given by $a_{1}, \ldots, a_{n}$.

\begin{definition}
An $Sp(2n,\mathbb{C})$-Higgs bundle is a pair $(E,\Phi)$ for $E$ a rank $2n$ vector bundle  with a symplectic form $\omega(~,~)$, and the Higgs field $\Phi\in H^{0}(\Sigma, \text{End}(E)\otimes K)$ satisfying
$$\omega(\Phi v,w)=-\omega(v,\Phi w).$$ 
\end{definition}
The volume form $\omega^{n}$ trivialises the determinant bundle $\Lambda^{2n}E^{*}$. The characteristic polynomial $\text{det}(\eta-\Phi)$ defines a spectral curve $\pi: S\rightarrow \Sigma$ in  $X$ with equation
\begin{eqnarray}
 \eta^{2n}+a_{1}\eta^{2n-1}+\ldots+a_{n-1}\eta^{2}+a_{n}=0, \label{eqsim}
\end{eqnarray}
whose genus  is $g_{S}:=1+4n^{2}(g-1).$ The  curve $S$ has a natural involution $\sigma(\eta)=-\eta$ and thus one can define the quotient curve $\bar \pi:\overline{S}=S/\sigma \rightarrow \Sigma$, of which $S$ is a 2-fold cover
\[p: S\rightarrow \overline{S}.\]
Note that   the Norm map associated to $p$  satisfies $p^{*}\text{Nm}(x)=x+\sigma x$, and thus the Prym variety $\text{Prym}(S,\overline{S})$ is given by the line bundles  $M\in \text{Jac}(S)$ for which $\sigma^{*}M\cong M^{*}$.

As in the case of classical Higgs bundles,  the characteristic polynomial of a Higgs field $\Phi$ gives the Hitchin fibration
\begin{eqnarray}
 h~:~ \mathcal{M}_{Sp(2n,\mathbb{C})}\longrightarrow\mathcal{A}_{Sp(2n,\mathbb{C})}:=\bigoplus_{i=1}^{n}H^{0}(\Sigma,K^{2i}),
\end{eqnarray}
and one has the following:
\begin{framed}
The generic fibres $h^{-1}(a)$ of the Hitchin fibration for $Sp(2n,\mathbb{C})$-Higgs bundles is given by Prym varieties $\text{Prym}(S,\overline{S})$ where $S$ and its quotient $\bar S$ are the curves  defined by $a$ as above.
\end{framed}

The spectral data described above for an $Sp(2n,\mathbb{C})$-Higgs bundle $(E,\Phi)$ can be obtained by looking at the extra conditions needed on  $L\in \text{Jac}(S)$ associated to the  corresponding classical Higgs pair for which $\pi_*L=E$. In order to understand this, note that for $\mathcal{V}\subset S$ an open set, we have $\mathcal{V}\subset \pi^{-1}(\pi(\mathcal{V}))$ and hence a natural restriction map
$  H^0(\pi^{-1}(\pi(\mathcal{V})), L) \rightarrow H^0(\mathcal{V},L),$
which gives the evaluation map $ev: ~\pi^*\pi_*L\rightarrow L$. Multiplication by $\eta$  commutes with this linear map and so the action of $\pi^*\Phi$ on the dual of the vector bundle $\pi^*\pi_*L$ preserves a one-dimensional subspace. Hence $L^*$ is an eigenspace of $\pi^*\Phi^t$, with eigenvalue $\eta$. Equivalently, $L$ is the cokernel of $\pi^*\Phi-\eta$ acting on $\pi^*E\otimes \pi^{*}K^{*}$. By means of the Norm map for $\pi$, this correspondence can be seen on  the curve $S$ via the exact sequence  \begin{eqnarray}
 0\rightarrow L\otimes \pi^{*}K^{1-2n}\rightarrow \pi^{*}E\xrightarrow{\pi^{*}\Phi-\eta} \pi^{*}(E\otimes K^*)\xrightarrow{ev}L\otimes \pi^{*}K\rightarrow0,\label{sequence1}
\end{eqnarray} 
and its dualised sequence 
\begin{eqnarray}
 0\rightarrow L^{*}\otimes \pi^{*}K^{*}\rightarrow \pi^{*}(E^{*}\otimes K^{*})\rightarrow \pi^{*}E^{*}\rightarrow L^{*}\otimes \pi^{*}K^{2n-1}\rightarrow 0. \label{sequence1dual}
\end{eqnarray} 
In particular, from the relative duality theorem one has that
\begin{eqnarray}
 \pi_{*}(L)^{*}\cong \pi_{*}(K_{S}\otimes \pi^{*}K^{-1}\otimes L^{*}),
\end{eqnarray}
and thus $E^{*}$ is the direct image sheaf $\pi_{*}(L^{*}\otimes \pi^{*}K^{2n-1})$.

Given  an $Sp(2n,\mathbb{C})$-Higgs bundle $(E,\Phi)$, one has $\Phi^{t}=-\Phi$ and  an eigenspace $L$ of $\Phi$ with eigenvalue $\eta$ is transformed to $\sigma^{*}L$ for the eigenvalue $-\eta$. Moreover, since the line bundle $L$ is the cokernel of $\pi^{*}\Phi-\eta$ acting on $\pi^{*}(E\otimes K^{*})$, one can consider the corresponding exact sequences (\ref{sequence1}) and its dualised sequence, which identify $L^{*}$ with $L \otimes \pi^{*}K^{1-2n}$, or equivalently, $L^{2} = \pi^{*}K^{2n-1}.$ By choosing a square root $K^{1/2}$ one has a line bundle $M:=L\otimes \pi^{*}K^{-n+1/2}$ for which $\sigma^{*}M\cong M^{*}$, i.e., which is in the Prym variety  $\text{Prym}(S,\overline{S})$.

Conversely, an $Sp(2p,\mathbb{C})$-Higgs bundle can be recovered from a line bundle $M\in \text{Prym}(S,\overline{S})$, for $S$ a smooth curve with equation (\ref{eqsim}) and $\bar S$ its quotient curve. Indeed, by Bertini's theorem, such a smooth curve $S$ with equation (\ref{eqsim}) always exists. Letting $E:=\pi_{*}L$ for $L=M\otimes \pi^{*}K^{n-1/2}$, one has the exact sequences (\ref{sequence1}) and its dualised on the curve $S$.  Moreover, since $L^{2}\cong \pi^{*}K^{2n-1}$, there is an isomorphism $E\cong E^{*}$ which induces the symplectic structure on $E$. Hence,  the generic fibres of the corresponding Hitchin fibration can be identified with the Prym variety $\text{Prym}(S,\overline{S})$.

\subsection{$SO(2n+1,\mathbb{C})$-Higgs bundles  } 

We shall now consider the special orthogonal group $G_{c}=SO(2n+1,\mathbb{C})$ and the corresponding Higgs bundles. 
Following a similar analysis as in the previous case, one can see that for a generic matrix $A\in \mathfrak{so}(2n+1,\mathbb{C})$, its  distinct eigenvalues occur in $\pm \lambda_{i}$ pairs. Thus, the characteristic polynomial of $A$ must be of the form
\begin{eqnarray}\text{det}(x-A)=x(x^{2n}+a_{1}x^{2n-2}+\ldots +a_{n-1}x^{2}+a_{n}),\label{charpolso}\end{eqnarray}
where the coefficients $a_{1}, \ldots, a_{n}$ give a basis for the invariant polynomials  on  $\mathfrak{so}(2n+1,\mathbb{C})$. 
\begin{definition}
An $SO(2n+1,\mathbb{C})$-Higgs bundle is a pair $(E,\Phi)$ for $E$ a holomorphic vector bundle of rank $2n+1$ with  a non-degenerate symmetric bilinear form $(v,w)$, and  $\Phi$ a Higgs field  in $H^{0}(\Sigma,\text{End}_{0}(E)\otimes K)$ which satisfies
$$(\Phi v,w)=-(v,\Phi w).$$
\end{definition}
The moduli space $\mathcal{M}_{SO(2n+1,\mathbb{C})}$  has two connected components, characterised by a class $w_{2} \in H^{2}(\Sigma,\mathbb{Z}_{2}) \cong  \mathbb{Z}_{2}$, depending on whether $E$ has a lift to a spin bundle or not. The spectral curve induced by the characteristic polynomial in (\ref{charpolso}) is a reducible curve: an $SO(2n+1,\mathbb{C})$-Higgs field $\Phi$  always has a zero eigenvalue, and from \cite[Section 4.1]{N3} the zero eigenspace $E_{0}$ is given by $E_{0}\cong K^{-n}$.

From (\ref{charpolso}), the characteristic polynomial $\text{det}(\eta- \Phi)$  defines a component of the spectral curve, which for convenience we shall denote by $\pi:S\rightarrow \Sigma$, and  whose  equation is given by 
$\eta^{2n}+a_{1}\eta^{2n-2}+\ldots +a_{n-1}\eta^{2}+a_{n}=0,$
where $a_{i}\in H^{0}(\Sigma, K^{2i})$. This is a $2n$-fold cover of $\Sigma$, with genus $g_{S}= 1+4n^{2}(g-1).$ 
The Hitchin fibration in this case is given by the map
\begin{eqnarray}
 h~:~ \mathcal{M}_{SO(2n+1,\mathbb{C})}\longrightarrow\mathcal{A}_{SO(2n+1,\mathbb{C})}:=\bigoplus_{i=1}^{n}H^{0}(\Sigma,K^{2i}),
\end{eqnarray} which sends each pair $(E,\Phi)$ to the coefficients of $\text{det}(\eta-\Phi)$. As in the case of $Sp(2n,\mathbb{C})$, the  curve $S$ has an involution $\sigma$ which acts as $\sigma(\eta)=-\eta$. Thus, we may consider the quotient curve $\overline{S}=S/\sigma$ in the total space of $K^{2}$,  for which $S$ is a double cover
$p: S\rightarrow \overline{S}.$ In this case the regular fibres can be seen as follows:

\begin{framed}
The regular fibres $h^{-1}(a)$ of the $SO(2n+1,\mathbb{C})$ Hitchin fibration are given by Prym varieties $\text{Prym}(S,\bar S)$ together with a trivialization of each $M\in \text{Prym}(S,\bar S)$ over the zeros of $a_{n}$ defining $S$ as in (\ref{charpolso}).
\end{framed}

Following \cite{N3}, the symmetric bilinear form $(v,w)$ canonically defines a skew form $(\Phi v,w)$ on $E/E_{0}$ with values in $K$. 
Moreover, choosing a square root $K^{1/2}$ one can define
\[V=E/E_{0}\otimes K^{-1/2},\]
on which the corresponding skew form  is non-degenerate. The Higgs field $\Phi$ induces a transformation $\Phi'$ on $V$ which has characteristic polynomial 
\[\text{det }(x-\Phi')=x^{2n}+a_{1}x^{2n-2}+\ldots+ a_{n-1}x^{2}+a_{n}.\]
Note that this is exactly the case of $Sp(2n,\mathbb{C})$ described in Section \ref{subsec:sp}, and thus we may describe the above with a choice of a line bundle $M_{0}$ in the Prym variety $\text{Prym}(S,\overline{S})$. In particular, $S$ corresponds to the smooth spectral curve of an $Sp(2n,\mathbb{C})$-Higgs bundle.

When reconstructing the vector bundle $E$ with an $SO(2n+1,\mathbb{C})$ structure from an $Sp(2n,\mathbb{C})$-Higgs bundle $(V,\Phi')$ as in \cite[Section 4.3]{N3}, there is a mod 2 invariant associated to each zero of the 
coefficient $a_{n}$ of the characteristic polynomial $\text{det}(\eta-\Phi')$. This data comes from choosing a trivialisation of $M\in \text{Prym}(S,\overline{S})$ over the zeros of $a_{n}$, and defines a covering $P'$ of the Prym variety $\text{Prym}(S,\overline{S})$. The covering has two components corresponding to the spin and non-spin lifts of the vector bundle.  The identity component of $P'$, which corresponds to the spin case, is isomorphic to the dual of the symplectic Prym variety, and this is the generic fibre of the $SO(2n+1,\mathbb{C})$ Hitchin map.

\subsection{$SO(2n,\mathbb{C})$-Higgs bundles } 

Lastly, we consider $G_{c}=SO(2n,\mathbb{C})$. As in previous cases, the distinct eigenvalues  of a matrix  $A\in \mathfrak{so}(2n,\mathbb{C})$  occur in  pairs $\pm \lambda_{i}$, and thus the characteristic  polynomial of $A$ is of the form
$\text{det}(x-A)=x^{2n}+a_{1}x^{2n-2}+\ldots+a_{n-1}x^{2}+a_{n}.$
In this case the coefficient $a_{n}$ is the square of a polynomial $p_{n}$, the Pfaffian, of degree $n$. A basis for the invariant polynomials on the Lie algebra $\mathfrak{so}(2n,\mathbb{C})$ is 
$a_{1},a_{2},\ldots, a_{n-1}, p_{n},$
(the reader should refer, for example, to \cite{asla} and references therein for further details). 

\begin{definition}
An $SO(2n,\mathbb{C})$-Higgs bundle is a pair $(E,\Phi)$, for $E$ a holomorphic vector bundle of rank $2n$ with a non-degenerate symmetric bilinear form $(~,~)$, and    $\Phi\in H^{0}(\Sigma,\text{End}_{0}(E)\otimes K)$   satisfying
$(\Phi v,w)=-(v,\Phi w).$         
\end{definition}

Considering the characteristic polynomial $\text{det}(\eta-\Phi)$ of a Higgs bundle $(E,\Phi)$ one obtains  a $2n$-fold cover $\pi:S\rightarrow \Sigma$ whose equation is given by
\[\text{det}(\eta-\Phi)=\eta^{2n}+a_{1}\eta^{2n-2}+\ldots+a_{n-1}\eta^{2}+p_{n}^{2},\]
for $a_{i}\in H^{0}(\Sigma,K^{2i})$ and $p_{n}\in H^{0}(\Sigma, K^{n})$. Note that this curve has always singularities, which are given by   $\eta=0$.
The  curve $S$ has a natural involution $\sigma(\eta)=-\eta$, whose fixed points in this case are the singularities of $S$. 
The \textit{virtual} genus of $S$ can be obtained via the adjunction formula, giving
$g_{S}=1+4n^{2}(g-1).$

In order to define the spectral data, one may consider its non-singular model $\hat{\pi}:\hat{S}\rightarrow \Sigma$, whose genus is given by 
\begin{eqnarray}
 g_{\hat{S}}&=& g_{S}-\# \text{singularities} \nonumber\\
&=& 1+4n^{2}(g-1) -2n(g-1)\nonumber \\
&=& 1+2n(2n-1)(g-1).\nonumber
\end{eqnarray}
As the fixed points of $\sigma$ are double points, the involution extends to an involution $\hat{\sigma}$ on $\hat{S}$ which does not have fixed points.
Considering the associated basis of invariant polynomials for each Higgs field $\Phi$, one may define the Hitchin fibration
\begin{eqnarray}
 h~:~ \mathcal{M}_{SO(2n,\mathbb{C})}\longrightarrow\mathcal{A}_{SO(2n,\mathbb{C})}:=H^{0}(\Sigma,K^{n})\oplus \bigoplus_{i=1}^{n-1}H^{0}(\Sigma,K^{2i}).
\end{eqnarray}

In this case the line bundle associated to an $SO(2n,\mathbb{C})$-Higgs bundle is defined on the desingularisation $\hat{S}$ of $S$:
\begin{framed}
The smooth fibres $h^{-1}(a)$ of the $SO(2n,\mathbb{C})$ Hitchin fibration are given by $\text{Prym}(\hat{S},\hat{S}/\hat{\sigma})$, for $\hat S$ the desingularisation of the curve $S$ associated to the regular base point $a$.
\end{framed}

Starting with an $SO(2n,\mathbb{C})$-Higgs bundle, since $\hat{S}$ is smooth we obtain an eigenspace bundle 
$L \subset  \text{ker}(\eta-\Phi)$ inside the vector bundle $E$ pulled back to $\hat{S}$.
 In particular, this line bundle satisfies
 $\hat{\sigma}^{*}L\cong L^{*}\otimes (K_{\hat{S}}\otimes \pi^*K^{*})^{-1},$
 thus defining a point in $\text{Prym}(\hat{S},\hat{S}/\hat{\sigma})$ given by
 \begin{eqnarray}
  M:=L\otimes (K_{\hat{S}}\otimes \pi^* K^{*})^{1/2}.\nonumber
 \end{eqnarray}

Conversely,   a Higgs bundle $(E,\Phi)$ may be recovered from a curve $S$ with has equation
$\eta^{2n}+a_{1}\eta^{2n-2}+\ldots+a_{n-1}\eta^{2}+p_{n}^{2}=0,$
and a line bundle $L$ on its desingularisation $\hat{S}$. Note that given the sections
\[s=\eta^{2n}+a_{1}\eta^{2n-2}+\ldots+a_{n-1}\eta^{2}+p_{n}^{2}\]
for fixed $p_{n}$ with simple zeros, one has a linear system whose only base points are when $\eta=0$ and $p_{n}=0$. Hence, by Bertini's theorem the generic divisor of the linear system defined by the sections $s$ has those base points as its only singularities. Moreover, as $p_{n}$ is a section of $K^{n}$, in general there are $2n(g-1)$ singularities which are generically  ordinary double points. A generic divisor of the above linear system defines a curve $S$ which has an involution $\sigma(\eta)=-\eta$ whose only fixed points are the base points. 

The involution $\sigma$ induces an involution $\hat{\sigma}$ on the desingularisation $\hat{S}$ of $S$ which has no fixed points, and thus we may consider the quotient $\hat{S}/\hat{\sigma}$ and the corresponding Prym variety $\text{Prym}(\hat{S},\hat{S}/\hat{\sigma})$. Following a similar procedure as before, a line bundle $M\in\text{Prym}(\hat{S},\hat{S}/\hat{\sigma}) $ induces a Higgs bundle $(E,\Phi)$ where $E$ is the direct image sheaf of $L=M\otimes (K_{\hat{S}}\otimes \pi^*K^{*})^{-1/2}$.
It is thus the Prym variety of $\hat{S}$ which is a generic fibre of the corresponding Hitchin fibration.  

\begin{exercise}
Show that the genus $ g_{\hat S/\hat \sigma}$ of $\hat S/\hat \sigma$ is $n(2n-1)(g-1)$.
\end{exercise}

\nobreak

\chapter{Spectral data for $G$-Higgs bundles}\label{cap2}
 \hfill{\vbox{\hbox{\it But most of all a good example is }\hbox{\it a thing of beauty. It shines and }\hbox{\it convinces. It gives insight and}\hbox{\it understanding. It provides the}\hbox{\it bedrock of belief.}\hbox{\it\noindent\rule{6cm}{0.4pt}}}}
 
 \hfill{\vbox{\hbox{Sir Michael Atiyah}}}

\bigskip

 Higgs bundles for real forms were first studied by N. Hitchin in \cite{N1}, and the results for $SL(2,\mathbb{R})$ were generalised in \cite{N5}, where Hitchin studied the case of $G=SL(n,\mathbb{R})$. Using Higgs bundles he counted the number of connected components and, in the case of split real forms, he identified a component homeomorphic to $\mathbb{R}^{\text{dim} G(2g-2)}$ and which naturally contains a copy of a  Teichm\"uller space. 
The aim of this Lecture is to introduce principal Higgs bundles for real forms and their corresponding spectral data as studied in \cite{thesis} and further developed in \cite{non ab, classes}. We begin by reviewing definitions and properties related to real forms of Lie algebras and Lie groups (see e.g., \cite{sym, helga, lie, knapp, int}), and then define $G$-Higgs bundles for real forms $G$ of classical complex Lie groups $G_c$. Through the approach of \cite{N5}, we describe these Higgs bundles as the fixed points of  a certain involution on the moduli space of $G_{c}$-Higgs bundles.  In later sections we study $G$-Higgs bundles for  non-compact real forms $G$ and in each case give an overview of the corresponding spectral data when available.

 \section{$G$-Higgs bundles}

A theorem by Hitchin \cite{N2} and Simpson \cite{simpson88} gives the most important property of stable Higgs bundles on a compact Riemann surface $\Sigma$ of genus $g\geq2$: 

\begin{theorem}
 If a Higgs bundle $(E,\Phi)$ is stable and  $\text{deg} ~E = 0$, then there is a unique unitary connection $A$ on $E$, compatible with the holomorphic structure, such that
\begin{eqnarray}
 F_{A}+ [\Phi,\Phi^{*}]=0~\in \Omega^{1,1}(\Sigma, \text{End}~E), \label{2.1}
\end{eqnarray}
 where $F_{A}$ is the curvature of the connection.  
\end{theorem}
Equation (\ref{2.1}) together with the holomorphicity condition
$
d_{A}''\Phi=0 \label{hit2}
$
give the \textit{Hitchin equations}, where  $d_{A}''\Phi$ is the anti-holomorphic part of the covariant derivative of $\Phi$. Following \cite{N5}, the above equations can also be considered when $A$ is a connection on a principal $G$-bundle $P$, where $G$ is the compact real form of a complex Lie group $G_{c}$, and $a\rightarrow -a^{*}$ is the compact real structure on the Lie algebra. This motivates the study of Higgs bundles for real forms.

\subsection{Real forms}\label{split}

Let   $\mathfrak{g}_{c}$ be a complex Lie algebra with complex structure $i$, whose Lie group is $G_{c}$.  
\begin{definition}
 A \em{real form} of $\mathfrak{g}^{c}$ is a real Lie algebra which satisfies
$\mathfrak{g}^{c}=\mathfrak{g}\oplus i\mathfrak{g}.$
\end{definition}

Given a real form $\mathfrak{g}$ of $\mathfrak{g}^{c}$, an element  $Z\in \mathfrak{g}^{c}$ may be written as $Z=X+iY$ for $X,Y\in \mathfrak{g}$. The mapping
$ X+iY\mapsto X-iY$
is called the \textit{conjugation}  with respect to $\mathfrak{g}$. 

\begin{remark}
 Any real form $\mathfrak{g}$ of $\mathfrak{g}^{c}$ is given by the fixed points set of an antilinear involution $\tau$ on $\mathfrak{g}^{c}$. 
  In particular the conjugation with respect to $\mathfrak{g}$ satisfies these properties. 
\end{remark}
%

\begin{definition}
 A \text{real form} of a complex Lie group $G_{c}$ is  an antiholomorphic Lie group automorphism $ \tau: G_{c}\rightarrow G_{c} $ of order two, i.e., 
$ \tau^{2}=1.$\end{definition}

Every  $X\in \mathfrak{g}^{c}$ defines an endomorphism $\text{ad}X$ of $\mathfrak{g}^{c}$ given by
$\text{ad}X(Y)=[X,Y] $ for $ Y\in \mathfrak{g}^{c}.$
For $\text{Tr}$ the trace of a vector space endomorphism,  $B(X,Y)=\text{Tr}(\text{ad}X\text{ad}Y)$ is  a the bilinear form on $\mathfrak{g}^{c}\times \mathfrak{g}^{c}$  called the \textit{Killing form} of $\mathfrak{g}^{c}$. 

\begin{definition}
 A real Lie algebra $\mathfrak{g}$ is called {\em  compact } if the Killing form is negative definite on it. The corresponding Lie group $G$ is a compact Lie group.
\end{definition}

\begin{definition}\label{defisplit} Let $\mathfrak{g}$ be a real form of a complex simple Lie algebra $\mathfrak{g}^{c}$, given by the fixed points of an antilinear involution $\tau$. Then, 
   if there is a Cartan subalgebra invariant under $\tau$ on which the Killing form is negative definite, the real form $\mathfrak{g}$ is called a {\em compact real form}.
             Such a compact real form of $\mathfrak{g}^{c}$ corresponds to a compact real form $G$ of $G_{c}$; 
  if there is an invariant Cartan subalgebra on which the Killing form is positive definite,  the form is called a {\em  split (or normal) real form}. The corresponding Lie group $G$ is the split real form of $G_{c}$.   \label{compactsplit} \label{realformsdef}
\end{definition}

Any complex semisimple Lie algebra  $\mathfrak{g}^{c}$ has a compact and a split real form which are unique up to conjugation via $\text{Aut}_{\mathbb{C}} \mathfrak{g}^{c}$ (e.g., for $\mathfrak{sl(n,\mathbb{C}})$ these are $\mathfrak{su(n)}$  and $\mathfrak{sl}(n,\mathbb{R})$ respectively). 

\begin{remark} Recall that all Cartan subalgebras $\mathfrak{h}$ of a finite dimensional Lie algebra $\mathfrak{g}$ have the same dimension. The rank of $\mathfrak{g}$ is defined to be this dimension, and a real form $\mathfrak{g}$ of a complex Lie algebra $\mathfrak{g}^{c}$ is split if and only if the real rank of $\mathfrak{g}$ equals the complex rank of $\mathfrak{g}^{c}$. 
\end{remark}


An involution $\theta$ of a real semisimple Lie algebra $\mathfrak{g}$ such that the symmetric  bilinear form
$B_{\theta}(X,Y)=-B(X,\theta Y)$ 
is positive definite is called a \textit{Cartan involution}. Any real semisimple Lie algebra has a Cartan involution, and any two Cartan  involutions $\theta_{1},\theta_{2}$ of $\mathfrak{g}$ are conjugate via an automorphism of $\mathfrak{g}$, i.e., there is a map $\varphi$ in $\text{Aut} \mathfrak{g}$ such that $\varphi \theta_{1}\varphi^{-1}=\theta_{2}$. The decomposition of $\mathfrak{g}$ into eigenspaces of a Cartan involution $\theta$ is called the \textit{Cartan decomposition} of $\mathfrak{g}$. 

\begin{proposition}[\cite{knapp}]
 Let $\mathfrak{g}^{c}$ be a complex semisimple Lie algebra, and $\rho$ the conjugation with respect to a  compact real  form $\mathfrak{u}$ of $\mathfrak{g}^{c}$. Then, $\rho$ is a Cartan involution.
 \label{conjucomp}
\end{proposition}

\begin{proposition}[\cite{helga}] \label{helga} Any non-compact real form $\mathfrak{g}$ of a complex simple Lie algebra $\mathfrak{g}^{c}$ can be obtained from a pair $(\mathfrak{u},\theta)$, for $\mathfrak{u}$ the compact real form of $\mathfrak{g}^{c}$  and $\theta$ an involution on $\mathfrak{u}$. \label{realmethod}
\end{proposition}

 For completion, we shall recall here the construction of real forms from \cite{helga}. Let $\mathfrak{h}$ be the $+1$-eigenspaces of $\theta$ and $i\mathfrak{m}$ the $-1$-eigenspace of $\theta$ acting on $\mathfrak{u}$, thus having \begin{eqnarray}\mathfrak{u}=\mathfrak{h} \oplus i\mathfrak{m},\label{decu}\end{eqnarray}
Since $ \mathfrak{g}^{c}
=  \mathfrak{h}\oplus \mathfrak{m} \oplus i (\mathfrak{h} \oplus \mathfrak{m}), $ %
 there is a natural non-compact real form $\mathfrak{g}$ of $\mathfrak{g}^{c}$ given by 
\begin{eqnarray}\mathfrak{g}=\mathfrak{h} \oplus \mathfrak{m}.\label{decnon}\end{eqnarray}
Moreover, if a linear isomorphism $\theta_{0}$ induces the  decomposition as in (\ref{decnon}), then $\theta_{0}$ is a Cartan involution of  $\mathfrak{g}$ 
and   $\mathfrak{h}$ is the maximal compact subalgebra of $\mathfrak{g}$.

Following the notation of Proposition \ref{realmethod}, let $\rho$ be the antilinear involution defining the compact form $\mathfrak{u}$ of a complex simple Lie algebra $\mathfrak{g}^{c}$ whose decomposition via an involution $\theta$ is given by equation (\ref{decu}).
Moreover, let $\tau$ be an antilinear involution which defines  the corresponding non-compact real form $\mathfrak{g}=\mathfrak{h}\oplus \mathfrak{m}$ of $\mathfrak{g}^{c}$. Considering the action of the two antilinear involutions $\rho$ and $\tau$ on $\mathfrak{g}^{c}$, we may decompose the Lie algebra $\mathfrak{g}^{c}$ into eigenspaces 
\begin{eqnarray}
 \mathfrak{g}^{c}=\mathfrak{h}^{(+,+)}\oplus \mathfrak{m}^{(-,+)}\oplus (i\mathfrak{m})^{(+,-)}\oplus (i\mathfrak{h})^{(-,-)},\label{decg}
\end{eqnarray}
where the upper index $(\cdot,\cdot)$ represents the $\pm$-eigenvalue of $\rho$ and $\tau$ respectively. 
 From the decomposition (\ref{decg}), the involution $\theta$ on the compact real form $\mathfrak{u}$  giving a non-compact real form $\mathfrak{g}$ of $\mathfrak{g}^{c}$ can be seen as acting on $\mathfrak{g}^{c}$ as  $
\sigma :=\rho\tau.$ 
Moreover, this induces an involution on the corresponding Lie group
$\sigma :=G_{c}\rightarrow G_{c}.$

\begin{remark}\label{fixedsigma2}
The fixed point set $\mathfrak{g}^{\sigma}$ of $\sigma$  is given by
$\mathfrak{g}^{\sigma}=\mathfrak{h}\oplus i\mathfrak{h}, $
and thus it is the complexification of the maximal compact subalgebra $\mathfrak{h}$ of $\mathfrak{g}$. Equivalently, the anti-invariant set under the involution $\sigma$ is given by $\mathfrak{m}^{\mathbb{C}}$.
\end{remark}

\subsection{$G$-Higgs bundles through involutions }\label{secinvo}
As mentioned previously, non-abelian Hodge theory on the compact Riemann surface $\Sigma$ gives a correspondence between the moduli space of reductive representations of $\pi_{1}(\Sigma)$ in a complex Lie group $G_{c}$ and the moduli space of $G_{c}$-Higgs bundles. The anti-holomorphic operation of conjugating by a real form $\tau$ of $G_{c}$ in the moduli space of representations can be seen via this correspondence as a holomorphic involution $\Theta$ of the moduli space of $G_{c}$-Higgs bundles.

Following \cite{N2}, in order to obtain a $G$-Higgs bundle, for $A$ the connection which solves Hitchin equations (\ref{2.1}), one requires the flat $GL(n,\mathbb{C})$ connection
\begin{eqnarray}
 \nabla=\nabla_{A} +\Phi +\Phi^{*}
\end{eqnarray}
to have  holonomy in a non-compact real form $G$ of $GL(n,\mathbb{C})$, whose real structure is $\tau$ and Lie algebra is $\mathfrak{g}$. Equivalently, for a complex Lie group $G_{c}$ with non-compact real form $G$ and real structure $\tau$, one requires
\begin{eqnarray}
\nabla=\nabla_A+\Phi-\rho(\Phi)
\end{eqnarray}
to have  holonomy in  $G$, where $\rho$ is the compact real structure of $G_{c}$. Since $A$ has holonomy in the compact real form of $G_{c}$, we have $\rho(\nabla_A)=\nabla_A$.   Hence, requiring  $\nabla=\tau(\nabla)$ is equivalent to  $\nabla_{A}=\tau(\nabla_{A})$
and $\Phi-\rho(\Phi)=\tau(\Phi-\rho(\Phi)).$
In terms of $\sigma=\rho\tau$, these two equalities are given by $\sigma(\nabla_{A})=\nabla_{A}$ and 
$ \Phi-\rho(\Phi)=\tau(\Phi-\rho(\Phi))
                =\tau(\Phi)-\sigma(\Phi)                =\sigma(\rho(\Phi)-\Phi). 
$
Hence, $\nabla$ has holonomy in the real form $G$ if $\nabla_{A}$ is invariant under  $\sigma$, and $\Phi$ anti-invariant. 
In terms of a $G_{c}$-Higgs bundle $(P,\Phi)$, one has that for $\mathcal{U}$ and $\mathcal{V}$  two trivialising open sets  in the compact Riemann surface $\Sigma$, the involution $\sigma$ induces an action on the 
transition functions $g_{uv}: \mathcal{U}\cap \mathcal{V}\rightarrow G_{c}$ given by
$g_{uv}\mapsto \sigma(g_{uv} ),$
and on the Higgs field by sending
$\Phi\mapsto -\sigma(\Phi).$
                              
Concretely,   for $G$ a real form of a complex semisimple lie group $G_{c}$, we may construct $G$-Higgs bundles as follows. For $H$ the maximal compact subgroup of $G$, we have seen that the Cartan decomposition of $\mathfrak{g}$ is given by
$\mathfrak{g}=\mathfrak{h}\oplus \mathfrak{m},$
for $\mathfrak{h}$ the Lie algebra of $H$, and $\mathfrak{m}$ its orthogonal complement. 
 This induces the following decomposition of the Lie algebra $\mathfrak{g}^{c}$ of $G_{c}$ in terms of the eigenspaces of the corresponding involution $\sigma$ as defined before:
$\mathfrak{g}^{c}=\mathfrak{h}^{\mathbb{C}}\oplus \mathfrak{m}^{\mathbb{C}}.$
Note that the Lie algebras satisfy
$ [\mathfrak{h}, \mathfrak{h}]\subset\mathfrak{h}$,  $[\mathfrak{h,\mathfrak{m}}]\subset\mathfrak{m}$, and $[\mathfrak{m},\mathfrak{m}]\subset \mathfrak{h}\nonumber
$. Hence there is an induced isotropy representation given by
$\text{Ad}|_{H^{\mathbb{C}}}: H^{\mathbb{C}}\rightarrow GL(\mathfrak{m}^{\mathbb{C}}). $
Then,   Definition \ref{principalLie} generalises to the following (see e.g. \cite{Go2}):

\begin{definition}
 A \text{principal} $G${\em-Higgs bundle} is a pair $(P,\Phi)$ where
  $P$ is a holomorphic principal $H^{\mathbb{C}}$-bundle on $\Sigma$, and 
  $\Phi$ is a holomorphic section of $P\times_{Ad}\mathfrak{m}^{\mathbb{C}}\otimes K$.
\label{defrealhiggs}
\end{definition}

\begin{example}\label{compactexample}
 For a compact real form $G$, one has $G=H$ and $\mathfrak{m}=\{0\}$, and thus $\sigma$ is the identity and the Higgs field must vanish:  a $G$-Higgs bundle becomes a principal $G_{c}$- bundle.
\end{example}

In terms of involutions, following \cite{N5} and recalling the previous analysis leading to Remark \ref{fixedsigma2}, we have the following:

\begin{proposition}\label{invrelation}
Let $G$ be a real form of a complex semi-simple Lie group $G_{c}$, whose real structure is $\tau$. Then,  $G${\em -Higgs bundles} are given by the fixed points in $\mathcal{M}_{G_{c}}$ of the involution $\Theta_G$ acting by
\[\Theta_G: ~(P,\Phi)\mapsto (\sigma(P),-\sigma(\Phi)),\]
where $\sigma=\rho\tau$, for $\rho$ the compact real form of $G_{c}$.     

\end{proposition} 
Similarly to the case of $G_{c}$-Higgs bundles, there is a notion of stability, semi stability and polystability for $G$-Higgs bundles. Following \cite[Section 3]{brad} and \cite[Section 2.1]{brad1}, one can see that the polystability of a $G$-Higgs bundle for $G\subset GL(n,\mathbb{C})$ is equivalent to the polystability of the corresponding $GL(n,\mathbb{C})$-Higgs bundle. However, a $G$-Higgs bundle can be stable as a $G$-Higgs bundle but not as a $GL(n,\mathbb{C})$-Higgs bundle.  We shall denote by $\mathcal{M}_{G}$ the moduli space of polystable $G$-Higgs bundles on   $\Sigma$.

\begin{exercise}[(*)]
Considering the notion of ``strong real form'' from \cite{adam}, describe the corresponding Higgs bundles and give a definition of $\Theta_{G}$ for which one does not have the problem described in the above paragraph. \end{exercise}

  One should note that  
a fixed point of $\Theta_{G}$ in $\mathcal{M}_{G^c}$ gives a representation of $\pi_1(\Sigma)$ into the real form $G$ up to the equivalence of conjugation by the normalizer of $G$ in $G^c$. This may be bigger than $G$ itself, and thus two distinct classes in $\mathcal{M}_{G}$ could be isomorphic in $\mathcal{M}_{G^c}$ via a complex map. 
Hence, although there is a map from $\mathcal{M}_{G}$ to the fixed point subvarieties in $\mathcal{M}_{G^c}$, this might not be an embedding. The reader may refer to \cite{GP09} for the Hitchin-Kobayashi type correspondence for real forms.
 
 \begin{remark}
A description of the above phenomena in the case of rank 2 Higgs bundles is given in \cite[Section 10]{Lau}, where one can see how the $SL(2,\mathbb{R})$-Higgs bundles which have different topological invariants lie in the same connected component as $SL(2,\mathbb{C})$-Higgs bundles. 
\end{remark}
\begin{remark}
The point of view of Proposition \ref{invrelation}, which is considered throughout \cite{thesis}, fits into a more global picture where $\Theta_G$ is one of three natural involutions acting on the moduli space of Higgs bundles \cite{branes, real}, giving three families of $(B,A,A)$, $(A,B,A)$ and $(A,A,B)$ branes in $\mathcal{M}_{G_c}$ as the fixed point sets. One should note that the  fixed point sets of these involutions are of great importance when studying the relation of Langlands duality with Higgs bundles, as initiated in \cite{Tamas1, Kap} and \cite{N3}.
\end{remark}

\section{Spectral data for $G$-Higgs bundles}

As mentioned in the first Lecture, the moduli spaces $\mathcal{M}_{G_{c}}$ have a natural symplectic structure, which we  denoted by $\omega$. Moreover, following \cite{N2}, the involutions $\Theta_{G}$ send $\omega\mapsto -\omega$. Thus, at a smooth point, the fixed point set must be Lagrangian and so the expected dimension of $\mathcal{M}_{G}$ is half the dimension of $\mathcal{M}_{G_{c}}$. 
In order to describe the spectral data for real $G$-Higgs bundles, one   considers the moduli space $\mathcal{M}_{G}$ sitting inside $\mathcal{M}_{G_{c}}$  as fixed points of $\Theta_{G}$ in the Hitchin base $\mathcal{A}_{G_c}$ and the corresponding preserved fibres. 

 \begin{rem}\label{traces}
 Let $\mathfrak{g}^{c}$ be one of the classical Lie algebras $\mathfrak{sl}(n,\mathbb{C})$, $\mathfrak{so}(2n+1,\mathbb{C})$,  and $\mathfrak{sp}(2n,\mathbb{C})$. Then, for $\pi: \mathfrak{g}^{c}\rightarrow \mathfrak{gl}(V)$ a representation of $\mathfrak{g}^{c}$, the ring of invariant polynomials of $\mathfrak{g}^{c}$ is generated by $\text{Tr}(\pi(X)^{i})$, for $i\in \mathbb{N}$ and $X\in \mathfrak{g}^{c}$.
\end{rem}

By considering Cartan's classification of classical Lie algebras, we shall now describe $G$-Higgs bundles and their spectral data 
for   non-compact real forms of a classical complex  Lie algebra $\mathfrak{g}^{c}$.  For $I_{n}$ the unit matrix of order $n$, we denote by $I_{p,q},~ J_{n}$  and $K_{p,q}$ the matrices
\begin{small}
\begin{eqnarray}I_{p,q}=\left(
\begin{array}
 {cc}
-I_{p}&0\\
0&I_{q}
\end{array}\right)
,~~
J_{n}=\left(
\begin{array}
 {cc}
0&I_{n}\\
-I_{n}&0
\end{array}
\right)
,~~  
K_{p,q}=\left(
\begin{array}
 {cccc}
-I_{p}&0&0&0\\
0&I_{q}&0&0\\
0&0&-I_{p}&0\\
0&0&0&I_{q}
\end{array}
\right).\label{MatricesIJK}\end{eqnarray}
\end{small}
Following Proposition \ref{helga}, we shall study each complex Lie algebra $\mathfrak{g}_c$ and compact form $\mathfrak{u}$ with different involutions $\theta$ which give decompositions $\mathfrak{u}=\mathfrak{h}\oplus i\mathfrak{m}$. Then the corresponding natural non-compact real form is $\mathfrak{g}=\mathfrak{h}\oplus\mathfrak{m}$, and to make sense of Proposition \ref{invrelation} we  consider
the following Lie algebras, Lie groups, real forms, and holomorphic and anti holomorphic involutions:
\begin{table}[h]
\begin{center}
\begin{small}

\begin{tabular}{c|c|c|c|c|c}
  $\mathfrak{g}_{c}$ & Lie group $G_{c}$ & Split form &Compact  form $\mathfrak{u}$& anti-involution $\rho$ fixing $\mathfrak{u}$ & dim $\mathfrak{u}$  \\
\hline
$\mathfrak{a}_{n}$  & $SL(n,\mathbb{C})$&$\mathfrak{sl}(n,\mathbb{R})$	&$\mathfrak{su}(n)$  &$\rho(X)=-\overline{X}^{t}$& $n(n+1)$    \\
$\mathfrak{b}_{n}$ & $SO(2n+1,\mathbb{C})$&$\mathfrak{so}(n,n+1)$		&$\mathfrak{so}(2n+1)$ &$\rho(X)=\overline{X}$& $n(2n+1)$   \\
$\mathfrak{c}_{n}$  & $Sp(2n,\mathbb{C})$	&$\mathfrak{sp}(2n,\mathbb{R})$	&$\mathfrak{sp}(n)$    & $\rho(X)=J_{n}\overline{X}J_{n}^{-1}$&$n(2n+1)$   \\
$\mathfrak{d}_{n}$  & $SO(2n,\mathbb{C})$	&$\mathfrak{so}(n,n)$	&$\mathfrak{so}(2n)$ &  $\rho(X)=\overline{X}$& $n(2n-1)$   
\label{compact table}
 \end{tabular}
 \end{small}

\caption{Compact forms $\mathfrak{u}$ of classical Lie algebras}
\end{center}
\end{table}

\vspace{-0.3 in }
\begin{table}[h]
\begin{center}
\begin{small}
\begin{tabular}{c|c|c|c|c}
$\mathfrak{g}_{c}$  &  Real form  $G$&  $\tau$ fixing $G$ & Involution $\theta$ on $\mathfrak{u}$ & $\sigma=\rho\tau$  \\
\hline
$\mathfrak{a}_{n}$ 	&$SL(n,\mathbb{R})$ &$\rho(X)=-\overline{X}^{t}$ &$\theta(X)=\overline{X}.$&$\sigma(X) =  -X^{t} $ \\
  &$SU^{*}(2m)$&$\tau(X)=J_{m}\overline{X}J_{m}^{-1}$&$\theta(X)=J_{m}\overline{X}J_{m}^{-1}$&$\sigma(X) = -J_{m}X^{t}J_{m}^{-1}$\\
  &$SU(p,q)$& $\tau(X)=-I_{p,q}\overline{X}^{t}I_{p,q}$ &$\theta(X)=I_{p,q}X I_{p,q}$&$\sigma(X) = I_{p,q}X I_{p,q}$\\ \hline
  $\mathfrak{b}_{n}$     &$SO(p,q)$ &$ \tau(X)=I_{p,q}\overline{X}I_{p,q}.$&$\theta(X)=I_{p,q}X I_{p,q}$&$ \sigma(X) = I_{p,q}X I_{p,q} $\\\hline
 $\mathfrak{c}_{n}$ 	      &$Sp(2n,\mathbb{R})$&$\tau(X)=\overline{X}$&$\theta(X)=\overline{X} $& $\sigma(X) = J_{n}X J_{n}^{-1} $\\
  & $Sp(2p,2q)$ & $\tau(X)=-K_{p,q}X^{*}K_{p,q}.$&$\theta(X)=K_{p,q} X K_{p,q},$& $\sigma(X) =K_{p,q} X K_{p,q}$\\
\hline
$\mathfrak{d}_{n}$ 	&$SO(p,q)$ &$ \tau(X)=I_{p,q}\overline{X}I_{p,q}.$&$\theta(X)=I_{p,q}X I_{p,q}$&$ \sigma(X) = I_{p,q}X I_{p,q} $\\
 &$SO^{*}(2m)$  &$
 \tau(X)=J_{m}\overline{X}J_{m}^{-1}.
$&$\theta(X)=J_{m}\overline{X}J_{m}^{-1}$&  $\sigma(X) = J_{m}X J_{m}^{-1}$
\label{compact table2}
 \end{tabular}
 \end{small}

\caption{Non-compact forms $G$ of classical Lie algebras $G_c$}\label{table1}\label{invotable}
\end{center}
\end{table}

In the case of split real forms, following the methods of \cite{N5} one obtains a description of real Higgs bundles which we shall use in subsequent sections:

\begin{theorem}[\cite{thesis}] For $G$ the split real form of $G_c$, the fixed points of $\Theta_G$ in the smooth fibres of the Hitchin fibration for $G_c$-Higgs bundles are given by points of order two.\label{splitteo}
\end{theorem}

  \subsection{$G=SL(n,\mathbb{R})$-Higgs bundles } Higgs bundles for $SL(n,\mathbb{R})$ were first considered in \cite{N5}, where Hitchin studied a copy of Teichm\"uller space inside the moduli space of Higgs bundles for split real forms.  
 Following Definition \ref{defrealhiggs}, an $SL(n,\mathbb{R})$-Higgs bundle is a pair $(E,\Phi)$ where $E$ is a rank $n$ orthogonal vector bundle and $\Phi:E\rightarrow E\otimes K$ is a symmetric and traceless holomorphic map. 
%


\begin{proposition}\label{invosln}
   $SL(n,\mathbb{R})$-Higgs bundles are given by the fixed points of $$\Theta_{SL(n,\mathbb{R})}:(E,\Phi)\mapsto (E^{*},\Phi^{t})$$  in $\mathcal{M}_{SL(n,\mathbb{C})}$corresponding to automorphisms $f:E\rightarrow E^{*}$  giving  a symmetric  form on $E$.
\end{proposition}

\begin{exercise}
Find the decomposition of $\mathfrak{u}=\mathfrak{su}(n)$ induced by the corresponding $\theta$ in Table \ref{invotable}, and use this to deduce Proposition \ref{invosln}.
\end{exercise}

Recalling that the trace is invariant under transposition, one has that 
the ring of invariant polynomials of $\mathfrak{g}^{c}=\mathfrak{sl}(n,\mathbb{C})$ is acted on trivially by the involution $-\sigma$, and thus the Hitchin base is preserved by $\Theta_{SL(n,\mathbb{R})}$. In order to find the spectral data for $SL(n, \mathbb{R})$-Higgs bundles, following Theorem \ref{splitteo} we look at elements of order two in the fibres of the Hitchin fibration for $SL(n,\mathbb{C})$-Higgs bundles:

\begin{framed}
Over a smooth point in the Hitchin base $\mathcal{A}_{SL(n,\mathbb{C})}$, Higgs bundles for $SL(n,\mathbb{R})$ correspond to line bundles $L\in \text{Prym}(S,\Sigma)$ such that $L^2\cong \mathcal{O}$.\end{framed}

In the case of $n=2$, the $SL(2,\mathbb{C})$-spectral curve $S$ given as in  (\ref{spectralsl}) has a natural involution $\sigma: \eta\mapsto -\eta$ and $\text{Prym}(S,\Sigma)=\{L\in \text{Jac}(S) :  \sigma^* L\cong L^*\}$. Hence, points in the smooth fibres corresponding to $SL(2,\mathbb{R})$-Higgs bundles are given by line bundles $L\in \text{Jac}(S)$ such that $\sigma^*L\cong L$.\begin{exercise}\label{exsl}
Let $L\in \text{Prym}(S,\Sigma)$ be a line bundle of order two. Then, its direct image is a rank 2 bundle on $\Sigma$ which decomposes into the sum of two line bundles $V\oplus V^*$. How can the Lefschetz fixed point formula from \cite{at} be used to relate the degree of $V$  and the action of $\sigma$ on $L$ in the sprit of \cite{umm}? 
\end{exercise}
The topological invariant associated to $SL(n,\mathbb{R})$-Higgs bundles is the characteristic class $\omega_2\in \mathbb{Z}_2$ which is the obstruction to lifting the orthogonal bundle to a spin bundle, and its study   was carried through in \cite[Section 5]{classes}.

\begin{exercise} For $n=2$, use the approach of  \cite[Proposition 3]{classes} to relate  $\omega_2$ to the invariants in Exercise \ref{exsl}.
\end{exercise}

The spectral data of $SL(n,\mathbb{R})$-Higgs bundles gives a finite cover of the smooth locus of the Hitchin fibration. For $n=2$, an explicit description of the monodromy action whose orbits are the connected components of $\mathcal{M}_{SL(2,\mathbb{R})}$ is given in \cite{Lau}.

\begin{exercise}[(*)] How can the methods in \cite{Lau} be extended to study monodromy for $SL(n,\mathbb{R})$-Higgs bundles for $n\geq 3$?
\end{exercise}


\subsection{$SU^{*}(2m)$-Higgs bundles}

The group $SU^*(2m)$  is the subgroup of $SL(2m,\mathbb{C})$ which commutes with an antilinear automorphism $J$ of $\mathbb{C}^{2m}$ such that $J^2=-1$.   At the level of the Lie algebras we have that the involution $\theta$ decomposes $\mathfrak{u}=\mathfrak{h}
\oplus i\mathfrak{m}$ where
$\mathfrak{h} = \mathfrak{sp}(m).$
The induced non-compact real form $\mathfrak{g}=\mathfrak{h}
\oplus  \mathfrak{m}$ is
\[\mathfrak{g}=\mathfrak{su}^{*}(2m)=\left\{
\left(
\begin{array}{cc}
 Z_{1}&Z_{2}\\
-\overline{Z}_{2}&\overline{Z}_{1}
\end{array}
\right)~\left|
\begin{array}{c}
Z_{1}, Z_{2}~ m\times m \text{~ ~ complex~matrices, } \\
\text{Tr}Z_{1}+\text{Tr}\overline{Z}_{1}=0
\end{array}\right.
\right\}.\]

\begin{definition} An $SU^{*}(2m)$ Higgs bundle on $\Sigma$ is given by
\[(E,\Phi)~ \text{for}~ \left\{ \begin{array}{c}
                                E \text{ ~ a~rank}~ 2m ~\text{ ~vector ~bundle ~with ~a ~symplectic ~form~}\omega \\
\Phi \in H^{0}(\Sigma, \text{End}(E)\otimes K)~\text{~ traceless~ and~symmetric~with~respect~to~} \omega.
                               \end{array}\right\}\]
\end{definition}

%

\begin{proposition}
Isomorphism classes of $SU^{*}(2m)$-Higgs bundles are given by fixed points of the involution $$\Theta_{SU^*}:(E,\Phi)\mapsto (E^{*},\Phi^{t})$$ on  $SL(2m,\mathbb{C})$-Higgs bundles  corresponding to pairs which have an automorphism   $f:E\rightarrow E^{*}$ endowing it with a symplectic structure, and which trivialises its determinant bundle. 
\end{proposition}

As the trace is invariant under conjugation and transposition, one has that the involution  $-\sigma(X) = J_{m}X^{t}J_{m}^{-1}$ acts trivially on the ring of invariant polynomials of $\mathfrak{sl}(2m,\mathbb{C})$, and thus preserves the Hitchin base. 
The spectral data associated to $SU^*(2m)$-Higgs bundles $(E,\Phi)$ was studied in \cite[Section 3]{non ab}, and we shall describe here its main features.  

The characteristic polynomial of an $SU^*(2m)$-Higgs bundle $(E,\Phi)$ can be seen to be the square of a Pfaffian, $\text{det}(\eta- \Phi)=p(\eta)^2$ and thus all fixed points of $\Theta_{SU^*}$ lie over singular points of the $SL(2m,\mathbb{C})$ Hitchin fibration. 
With a slight abuse of notation, we denote by $S$ the spectral curve  in the total space of $K$ defined by $$p(\eta)=\eta^m+a_2\eta^{m-2}+\dots+a_m=0$$ where the coefficients $a_i\in H^0(\Sigma, K^i)$.  It is a ramified $m$-fold cover of $\Sigma$ whose ramification points are the zeros of $a_m$. As in the case of complex groups, we interpret $p(\eta)=0$ as the vanishing of a section of $\pi^*K^m$ over the total space of the canonical bundle $\pi:K\rightarrow \Sigma$, where $\eta$ is the tautological section of $\pi^*K$, and Bertini's theorem assures us that for generic $a_i$ the curve is nonsingular.

\begin{exercise}
 What is the genus $g_S$ of $S$?
\end{exercise}

 On the spectral curve $S$,  the cokernel of $(\eta-\Phi)$ is a rank two holomorphic vector bundle $V$ on S. Then, following  \cite{bobi} (and using $p(\Phi)=0$ instead of the Cayley-Hamilton theorem), we can identify $E$ with the direct image $\pi_*V$ and  $\Phi$  as the direct image of $\eta:V\rightarrow V\otimes \pi^*K$. 
 
 \begin{exercise}
Given a proper $\Phi$-subbundle $F\subset E$, the characteristic polynomial of $\Phi|_{F}$ divides the one of $\Phi$. Use this together with Grothendieck-Riemann-Roch to show that semi-stability of $V$ implies semi-stability of $(E,\Phi)$.
\end{exercise}
 
Following \cite[Proposition 1]{non ab}  one obtains a description of the spectral data:
 
 \begin{framed}
The fixed point set of $\Theta_{SU^*(2m)}$ in a smooth fibre   of the $SL(2m,\mathbb{C})$-Hitchin fibration is  the moduli space of semi-stable rank 2 vector bundles on $S$ with fixed determinant $\pi^*K^{m-1}$.
\end{framed}

\begin{exercise}
Follow the approach of $SL(n,\mathbb{C})$-Higgs bundles to note that by fixing the determinant of $V$ one obtains a trivialization of the determinant of $\pi_*V$ on $\Sigma$.
\end{exercise}

\subsection{$SU(p,q)$-Higgs bundles}

\begin{definition}
 An $SU(p,q)$-Higgs bundle over $\Sigma$ is a pair $(E,\Phi)$ where
 $E=W_1\oplus W_2$ for $W_1,W_2$ vector bundles over $\Sigma$ of rank $p$ and $q$ such that $\Lambda^{p}W_1\cong \Lambda^{q}W_2^{*}$, and the Higgs field $\Phi$ is given by
$        \Phi=\left( \begin{array}
          {cc} 0&\beta\\
\gamma&0
         \end{array}\right), \label{hig}
$ for $\beta:W_2\rightarrow W_1\otimes K$ and $\gamma:W_1\rightarrow W_2 \otimes K$.  
\end{definition}
\begin{exercise}
Find the decomposition   $\mathfrak{u}=\mathfrak{h}\oplus i\mathfrak{m}$  via the action of $\theta$ in Table \ref{compact table} and deduce that $\theta \rho $ is the anti-holomorphic involution fixing the non-compact real form $\mathfrak{u}(p,q)$. \end{exercise}

\begin{proposition}
    $SU(p,q)$-Higgs bundles are  fixed points of  $\Theta_{SU(p,q)}:(E,\Phi)\mapsto (E,-\Phi)$ on $SL(p+q,\mathbb{C})$-Higgs bundles corresponding to bundles $E$ which have an automorphism conjugate to $I_{p,q}$ sending $\Phi$ to $-\Phi$, and  whose $\pm 1$ eigenspaces have dimensions $p$ and $q$. 
\end{proposition}

The involution $-\sigma$ acts trivially on the polynomials of even degree. 
%
Whilst the spectral data is not known for $p\neq q$, in the case of $p=q$ it has been described in \cite[Chapter 6]{thesis} and \cite{umm} by looking at $U(p,p)$-Higgs bundles $(W_1\oplus W_2, \Phi)$, which when satisfying $\Lambda^p W_1\cong \Lambda^q W_2^*$ correspond to $SU(p,p)$-Higgs bundles. In this case, the characteristic polynomial   defines a spectral curve $\pi:S\rightarrow \Sigma$ through the equation
$\text{det}(\eta-\Phi)= \eta^{2p}+a_{2}\eta^{2p-2}+\ldots+a_{2p-2}\eta^2+a_{2p}=0,$
where $\eta$ is the tautological section of $\pi^*K$ and $a_i\in H^0(\Sigma,K^i)$. This is a $2p$-fold cover of $\Sigma$, ramified over the $4p(g-1)$ zeros of $a_{2p}$,  and has a natural involution $\eta\mapsto -\eta$ which has as fixed points the ramification points of the cover, and which by abuse of notation, we shall call $\sigma$.

The  involution $\sigma$ plays an important role when constructing the spectral data as described in \cite[Proposition 2.1]{umm}. A line bundle $L$ on $S$ which defines a classical Higgs bundle induces a $U(p,p)$-Higgs bundle if and only if $\sigma^*L\cong L $. In this case, at a fixed point $x \in S$ of the involution, there is a linear action of $\sigma$ on the fibre $L_x$  given by scalar multiplication of $\pm 1$.
This description of the spectral data  can be then seen in terms of Jacobians through \cite[Theorem 3.2]{umm}:

\begin{framed}
    The fixed point set of $\Theta_{U(p,p)}$ in a smooth fibre of the classical Hitchin fibration can be seen in terms of  pull backs of $\text{Jac}(S/\sigma)$.  \end{framed}

The topological invariants associated to $U(p,p)$-Higgs bundles $(W_1\oplus W_2, \Phi)$ are the degrees $\text{deg}(W_1)$ and   $\text{deg}(W_2)$, and  can be seen in terms of the degree of the line bundle $L$ on $S$ and the number of ramification points of $S$ over which the linear action of $\sigma$ on the fibre of $L$ is $-1$.

\begin{exercise} Use the Lefschetz fixed point formula in \cite[Theorem 4.12]{at}  to see that the parity of the degree of $L$ and the number of points over which $\sigma$ acts as $-1$ needs to be the same. \label{exumm}
\end{exercise}

\begin{exercise}Following \cite{brad}, the Toledo invariant $\tau(\text{deg}(W_1),\text{deg}(W_2))$ associated to $U(p,p)$-Higgs bundles is defined as $\tau(\text{deg}(W_1),\text{deg}(W_2)):=\text{deg}(W_1)-\text{deg}(W_2)$. Use Exercise \ref{exumm} to express the invariant in terms of fixed points of $\sigma$ and obtain natural bounds.  
\end{exercise}

In the case of $SU(p,p)$-Higgs bundles, for maximal Toledo invariant   the fixed point set of   
  $\Theta_{SU(p,p)}$ in a smooth fibre  of the $SL(2p,\mathbb{C})$-Hitchin fibration is given by a covering of $\text{Prym}(S/\sigma, \Sigma)$, the Prym variety of the quotient curve $S/\sigma$. For $SU(p,p+1)$-Higgs bundles, the Cameral data can be used to obtain the spectral data as seen in \cite{ana1}.

 \begin{exercise}[(*)]
 How can the methods from \cite{umm} together with the approach of \cite{Xia2} be used to obtain the spectral data for $SU(p,1)$-Higgs bundles?
 \end{exercise}

\begin{remark}
It is interesting to note that the spectral data for $SU(p,p)$-Higgs bundles from \cite{umm} is used to conjecture a dual space to $\mathcal{M}_{SU(p,p)}$ through Langlands duality in \cite{classes}.
\end{remark}

\subsection{$SO(p,q)$-Higgs bundles}\label{sec:sopq}
In this case,   if $p+q$ is even, $\mathfrak{g}$ is a split real form if and only if $p=q$;
  if $p+q$ is odd, $\mathfrak{g}$  is a split real form if and only if $p=q+1$. 
Whilst we shall give some details on the construction of $SO(p,q)$-Higgs bundles, for a more detailed description of the approach needed to understand groups with signature the reader should refer to the following section on $Sp(2p,2q)$-Higgs bundles.

\begin{proposition}
 $SO (p,q)$ Higgs bundles are fixed points of  $$\Theta_{SO(p,q)}:~(E,\Phi)\mapsto (E, -\Phi)$$
on the moduli space of $SO(p+q,\mathbb{C})$ corresponding to vector bundles $E$ which have an automorphism $f$ conjugate to $I_{p,q}$ sending $\Phi$ to $-\Phi$ and whose $\pm 1$ eigenspaces have dimensions $p$ and $q$.
\end{proposition}

\begin{exercise}
The associated involution $\theta$ from Table \ref{invotable} decomposes $\mathfrak{u}=   \mathfrak{h}\oplus i \mathfrak{m}$. Give an explicit description of $\mathfrak{m}$ and $\mathfrak{h}$ and of the real form $\mathfrak{g}=   \mathfrak{h}\oplus  \mathfrak{m}$.
\end{exercise}

The vector space $V$ associated to the standard representation of $\mathfrak{h}^{\mathbb{C}}$  can be decomposed into
$V=V_{p}\oplus V_{q}$, 
for $V_{p}$ and $V_{q}$ complex vector spaces of dimension $p$ and $q$ respectively, with orthogonal structures. 
 The maximal compact subalgebra of $\mathfrak{so}(p,q)$   is $\mathfrak{h}=\mathfrak{so}(p)\times \mathfrak{so}(q)$ and  the Cartan decomposition of   $\mathfrak{so}(p+q,\mathbb{C})$ is given by
$(\mathfrak{so}(p,\mathbb{C})\oplus \mathfrak{so}(q,\mathbb{C}))\oplus \mathfrak{m}^{\mathbb{C}},$ for
\[\mathfrak{m}=
\left\{
\left.
\left(
\begin{array}
 {cc}
0&X_{2}\\
X_{2}^{t}&0
\end{array}
\right)
\right|
X_{2} ~\text{  real~} p\times q \text{~ matrix}
\right\}.\]

\begin{definition}
 An $SO(p,q)$ Higgs bundle is a pair $(E,\Phi)$ where $E=V_{p}\oplus V_{q}$ for $V_{p}$ and $V_{q}$ complex vector spaces of dimension $p$ and $q$ respectively, with orthogonal structures, and the Higgs field is a section in $H^{0}(\Sigma, (\text{Hom}(V_{q},V_{p})\oplus \text{Hom}(V_{p},V_{q}))\otimes K)$  given by
\[\Phi=\left(\begin{array}{cc}
              0&\beta\\
\gamma&0
             \end{array}
\right)~\text{~ for~}~\gamma \equiv -\beta^\text{T},~ \text{and~} \beta^\text{T} \text{the orthogonal transpose of }\beta.\] 
\end{definition}

Since the ring of invariant polynomials of $\mathfrak{g}^{c}=\mathfrak{so}(2m+1,\mathbb{C})$ is generated by $\text{Tr}(X^{i})$ for $X\in \mathfrak{g}^{c}$, for $p+q=2m+1$ one has
that the induced action of the involution $\Theta_{SO(p,q)}$ is trivial on the ring of invariant polynomials of the Lie algebra $\mathfrak{so}(2m+1,\mathbb{C})$, i.e., when $p$ and $q$ have different parity. 

\begin{exercise}
In the case of $\mathfrak{so}(2m,\mathbb{C})$, for $2m=p+q$, the ring of invariant polynomials is generated by $\text{Tr}(X^{i})$ for $X\in \mathfrak{g}^{c}$ and $i<2m$, together with the Pfaffian $p_{m}$, which is of degree $m$. Under which conditions on $p$ and $q$ is the induced action of $\Theta_{SO(p,q)}$ trivial on the ring of invariant polynomials?
\end{exercise}

%

\begin{framed}
	 The spectral data for $SO(p,q)$-Higgs bundles when $p=q$ or $p=q+1$ can be seen through Theorem \ref{splitteo} from \cite{thesis}  as points of order two in the smooth fibres of the $SO(p+q,\mathbb{C})$-Hitchin fibration. 
	 \end{framed}
	 
	 In both cases a key ingredient is the double cover $p:S\rightarrow S/\sigma$ given by the spectral curve (the desingularised curve in the case of $SO(2n,\mathbb{C})$) over the quotient curve, which through $K$-theoretic methods allows one to express the topological invariants involved in terms of the action of $\sigma$ \cite{sonn}.
	 

\subsection{$G=SO^{*}(2m)$} 

 
   The action of  $\theta$ in Table \ref{compact table} decomposes  $\mathfrak{u}=\mathfrak{h}\oplus i\mathfrak{m}$ for
$\mathfrak{h}= \mathfrak{u}(m)\cong \mathfrak{so}(2m)\cap \mathfrak{sp}(m),
$ and \begin{eqnarray}
i\mathfrak{m}&=& 
\left\{
\left.
\left(
\begin{array}
 {cc}
X_{1}&X_{2}\\
X_{2}&-X_{1}
\end{array}
\right)
\right|
X_{1},~X_{2} \in \mathfrak{so}(m)
\right\},
\end{eqnarray}
and the induced non-compact real form $\mathfrak{g}=\mathfrak{h}\oplus \mathfrak{m}$ is
\[\mathfrak{g}=\mathfrak{so}^{*}(2m)=\left\{
\left(
\begin{array}{cc}
 Z_{1}&Z_{2}\\
-\overline{Z}_{2}&\overline{Z}_{1}
\end{array}
\right)~\left|
\begin{array}{c}
Z_{1}, Z_{2}~ m\times m \text{~ ~ complex~matrices } \\
Z_{1} ~\text{  skew~symmetric,~} Z_{2}~\text{Hermitian}
\end{array}\right.
\right\}.\]
 

\begin{proposition}
 $SO^{*}(2m)$-Higgs bundles are   fixed points of the involution $$\Theta_{SO^{*}(2m)}:~(E,\Phi)\mapsto (E, -\Phi)$$ on the moduli space of $SO(2m,\mathbb{C})$-Higgs bundles corresponding to vector bundles $E$ which have an orthogonal automorphism  $f$ conjugate to $J_{m}$, sending $\Phi$ to $-\Phi$ and which squares to $-1$, equipping $E$ with a symplectic structure.  
\end{proposition}

The vector space associated to the standard representation of $\mathfrak{h}^{\mathbb{C}}$ has an orthogonal and symplectic structure $J$. Since $J^{-1}=J^{t}$ and $J^{2}=-1$, the vector space  may be expressed in terms of the $\pm i$ eigenspaces of $J$ as
$V\oplus V^{*},$
for $V$ a rank $m$ vector space. Thus, we have the following definition:

\begin{definition}
 An   $SO^{*}(2m)$\text{-Higgs bundle} is given by a pair $(E,\Phi)$ where $E=V\oplus V^{*}$ for $V$ a rank $m$ holomorphic vector bundle, and where the Higgs field $\Phi$ is given by
\[\Phi=\left(
\begin{array}
{cc}
0&\beta\\
\gamma&0 
\end{array}
\right)~\text{~ ~for~}~ \left\{
\begin{array}
 {l}
\gamma:~V\rightarrow V^{*}\otimes K  \text{~ satisfying~}\gamma=-\gamma^{t} \\
\beta:~V^{*}\rightarrow V \otimes K \text{~ satisfying~}\beta=-\beta^{t} 
\end{array}
 \right. .\]
\end{definition}

As in the previous case, the involution induced action of $\Theta_{SO^*(2m)}$ is trivial on the ring of invariant polynomials of $\mathfrak{g}^{c}$. The spectral data for these Higgs bundles is studied in \cite{non ab}, and we shall give a short description bellow. 

 In order to understand the associated spectral data, one notes that $SO^*(2m)$-Higgs bundles $(E,\Phi)$ are in fact $SU^{*}(2m)$\text{-Higgs bundle} with extra conditions. Hence, one may define a natural $m$ cover of the Riemann surface $\pi:S\rightarrow\Sigma$ by taking   $$\sqrt{\text{char}(\Phi)}=\eta^m+a_2\eta^{2m-2}+\ldots+a_m,$$
 and a rank 2 vector bundle $V$ on $S$ whose direct image on $\Sigma$ gives $E$. Since in this case the equation of the spectral curve  only has even coefficients, there is a natural involution $\sigma:\eta\rightarrow-\eta$ and one may consider the induced action of $\sigma$ on $V$ and on its determinant bundle. In particular, from \cite[Proposition 2]{non ab}   the vector bundle $V$ gives an $SO^*(2m)$-Higgs bundle if and only if it is preserved by the involution and the induced action on it satisfies some conditions:

 \begin{framed}
The fixed point set of $\Theta_{SO^*(2m)}$ in a smooth fibre   of the $SO(2m,\mathbb{C})$-Hitchin fibration is given by the moduli space of semi-stable rank 2 vector bundles $V$ on $S$ with fixed determinant  $\pi^*K^{2m-1}$, whose induced action by $\sigma$ on the determinant bundle is trivial.
\end{framed}

\begin{exercise}
The relative duality theorem gives 
\[(\pi_*(V))^*\cong \pi_*(V^*\otimes K_S)\otimes K^*.\]
 Use this to see that in order to have $E\cong E^*$ through a skew form,  the action of $\sigma$ needs to be trivial on the determinant bundle of $V$ for $\pi_*V=E$.
 \end{exercise}
 \begin{exercise}[(*)]
Describe how the vector bundles of rank 2 from \cite{non ab} appear in the description of the connected components of $\mathcal{M}_{SO^*(2m)}$ in  \cite{brad4}.\end{exercise}

 \label{sec:real2}

\subsection{$G=Sp(2n,\mathbb{R})$-Higgs bundles}\label{sec:sp1}

In this section and the one which follows we consider the non-compact real forms of the complex Lie group $Sp(2n,\mathbb{C})$. For this, recall that the symplectic Lie algebra $\mathfrak{sp}(2n,\mathbb{C})$ is given by the set of $2n\times 2n$ complex matrices $X$ that satisfy $J_{n}X + X^{t}J_{n} = 0 $ or equivalently,
$X=- J^{-1}_{n}X^{t}J_{n}.$

Let $\mathfrak{u}$ be the compact real form $\mathfrak{u}=\mathfrak{sp}(n)$ and $\theta(X)=\overline{X} =J_{n}\overline{X}J_{n}^{-1}.$  The Lie algebra $\mathfrak{sp}(n)$ is given by the quaternionic skew-Hermitian matrices; that is, the set of $n\times n$ quaternionic matrices $X$ which satisfy
$X=-\overline{X}^{t}$.
The compact form  may be decomposed as $\mathfrak{u}=\mathfrak{h}\oplus i\mathfrak{m}$, for
$\mathfrak{h}= 
\mathfrak{u}(n)\cong \mathfrak{so}(2n)\cap \mathfrak{sp}(n),$
which leads to the split real form $\mathfrak{g}=\mathfrak{h}\oplus   \mathfrak{m} $
 defined by \[\mathfrak{g}=\mathfrak{sp}(2n,\mathbb{R})=\left\{
\left(
\begin{array}{cc}
 X_{1}&X_{2}\\
X_{3}&-X^{t}_{1}
\end{array}
\right)~\left|
\begin{array}{c}
X_{i}~\text{  real~} n\times n \text{~ ~matrices } \\
X_{2}, X_{3}~ \text{  symmetric}
\end{array}\right.
\right\}.\]


\begin{definition}
 An   $Sp(2n,\mathbb{R})$\text{-Higgs bundle} is given by a pair $(E,\Phi)$ where $E=V\oplus V^{*}$ for $V$ a rank $n$ holomorphic vector bundle, and for $\Phi$ the Higgs field given by
\[\Phi=\left(
\begin{array}
{cc}
0&\beta\\
\gamma&0 
\end{array}
\right)~
\text{~ ~for~}~ \left\{
\begin{array}
 {l}
\gamma:~V\rightarrow V^{*}\otimes K \text{~ satisfying~}\gamma=\gamma^{t} \\
\beta:~V^{*}\rightarrow V\otimes K \text{~ satisfying~}\beta=\beta^{t} 
\end{array}
 \right. .\]
 
\end{definition}

\begin{proposition}
  $Sp(2n,\mathbb{C})$ Higgs bundles, and  $Sp(2n,\mathbb{R})$-Higgs bundles are given by the fixed points of $$\Theta_{Sp(2n,\mathbb{R})}:~(E,\Phi)\mapsto (E, -\Phi)$$ on $Sp(2n,\mathbb{C})$-Higgs bundles 
corresponding to vector bundles $E$ which have a symplectic isomorphism sending $\Phi$ to $-\Phi$, and whose square is the identity, endowing $E$ with an orthogonal structure. 
\end{proposition}


The invariant polynomials of $\mathfrak{g}^{c}$ are of even degree, and hence the involution $-\sigma$  acts trivially on them, making $\Theta_{Sp(2n,\mathbb{R})}$ preserve the whole Hitchin base $\mathcal{A}_{Sp(2n,\mathbb{C})}$. In the case of rank 4 Higgs bundles, the spectral data was first consider in P. Gothen's thesis \cite{Go2}, and through Theorem \ref{splitteo} and the spectral data for complex Higgs bundles one has the following:
\begin{framed}
The fixed points of $\Theta_{Sp(2n,\mathbb{R})}$ in the smooth fibres of the $Sp(2n,\mathbb{C})$-Hitchin fibration are given by line bundles $L\in \text{Prym}(S, S/\sigma)$ such that $L^{2}\cong \mathcal{O}$.
\end{framed}

In particular, since $S$ is a ramified double cover of $S/\sigma$, one has that  $L\in \text{Prym}(S,S/\sigma)$ if and only if $\sigma^*L\cong L^*$. Hence, by considering points of order two one has that $\sigma^*L\cong L$ and thus there is a natural induced action of $\sigma$ on the line bundle $L$. The topological invariants associated to these Higgs bundles were studied in   \cite{classes} through the natural action of $\sigma$.

\begin{exercise}
Compare the calculations in \cite[Section 4.1]{Go2} which lead to Milnor-Wood type inequalities for $Sp(2n,\mathbb{R})$-Higgs bundles, with the inequalities one obtains by using the involution $\sigma$ as in \cite{classes}.
\end{exercise}

\begin{exercise}[(*)]
Express the invariants from \cite[Section 4.2.4-4.2.6]{Go2} in terms of different choices of the natural involution $\sigma$ on $S$ as well as in terms of the action of a second natural involution appearing in some situations on $S/\sigma$. 
\end{exercise}

 \subsection{$G=Sp(2p,2q)$-Higgs bundles }\label{sec:sp2}

The induced non-compact real form $\mathfrak{g}=\mathfrak{h}\oplus  \mathfrak{m} $ is
\[\mathfrak{sp}(2p,2q)=\left\{
\left(
\begin{array}{cccc}
 Z_{11}&Z_{12}&Z_{13}&Z_{14}\\
 \overline{Z}^{t}_{12}&Z_{22}&Z^{t}_{14}&Z_{24}\\
 -\overline{Z}_{13}& \overline{Z}_{14}& \overline{Z}_{11}& -\overline{Z}_{12}\\
 \overline{Z}^{t}_{14}& -\overline{Z}_{24}&-Z^{t}_{12}& \overline{Z}_{22}\\
\end{array}
\right)~\left|
\begin{array}{c}
 Z_{i,j}~ \text{~ ~ complex~matrices, } \\
Z_{11}, ~Z_{13} \text{~ order~ } p,\\
Z_{12}, ~Z_{14} ~p\times q~\text{  matrices},\\
Z_{11}, ~Z_{22} ~\text{  skew ~Hermitian},\\
Z_{13}, ~Z_{24}~\text{  symmetric}.
\end{array}\right.
\right\}.\] 
\begin{exercise} Show that $\mathfrak{m}^{\mathbb{C}}$ can be expressed as subset of certain off-diagonal matrices.
\end{exercise}

\begin{proposition}
  $Sp(2p,2q)$-Higgs bundles are given by the fixed points of $$\Theta_{Sp(2p,2q)}:~(E,\Phi)\mapsto (E, -\Phi^{\text{T}})$$
on the moduli space of $Sp(2p+2q,\mathbb{C})$-Higgs bundles corresponding to symplectic vector bundles $E$ which have an endomorphism $f:E\rightarrow E$ conjugate to $\tilde{K}_{p,q}$, sending $\Phi$ to the symplectic transpose $-\Phi^\text{T}$, and whose $\pm 1$ eigenspaces are of dimension $2p$ and $2q$. 
\end{proposition}

\begin{definition}
 An $Sp(2p,2q)$\text{-Higgs bundle} is given by a pair $(E,\Phi)$ for $E=V_{2p}\oplus V_{2q}$ is a direct sum of symplectic vector spaces of rank $2p$ and $2q$, and    
\[\Phi=
\left(\begin{array}{cc}
0& -\gamma^\text{T}\\
\gamma&0
        \end{array}\right)~\text{  for~}~\left\{
\begin{array}{cc}
\gamma&: V_{2p}\rightarrow V_{2q}\otimes K  \\
-\gamma^\text{T}&: V_{2q}\rightarrow V_{2p} \otimes K
\end{array} \right., ~\text{for}~ \gamma^\text{T} \text{the symplectic transpose of} ~\gamma.
\]
\end{definition}

As the trace is invariant under conjugation and transposition, the induced action of $\Theta_{Sp(2p,2q)}$ is trivial on the ring of invariant polynomials of $\mathfrak{g}^{c}=\mathfrak{sp}(2(p+q),\mathbb{C})$.
In the case of $p=q$, one can see that $Sp(2p,2p)$-Higgs bundles are a  particular case of $SU^*(2p)$-Higgs bundles, and thus one needs to understand which extra conditions to the spectral data for  $SU^*(2p)$-Higgs bundles needs to be added in order to have the Higgs bundles for the symplectic real form. 

From the previous section, when $p=q$ the corresponding spectral curve is a $2p$-fold cover of the Riemann surface $\Sigma$ whose equation is given by the square root of the characteristic polynomial of the Higgs field. Moreover, it has a natural involution $\sigma$ whose action determines the associated spectral data. More precisely, the following  is shown in \cite[Chapter 7]{thesis}, and reinterpreted in \cite[Section 5]{non ab}:

 \begin{framed}
The fixed point set of $\Theta_{Sp(2p,2p)}$ in a smooth fibre   of the $Sp(4p,\mathbb{C})$-Hitchin fibration is given by the moduli space of semi-stable rank 2 vector bundles $V$ on $S$ with fixed determinant  $\pi^*K^{2p-1}$, whose induced action by $\sigma$ on $\Lambda^2V$  is $-1$.
\end{framed}

Since the action on $\Lambda^2V$ is $-1$, the involution $\sigma:S\rightarrow S$ acts  with different  eigenvalues $\pm 1$ on the fibres of $V$ over the ramification points of $S$, and thus through \cite{AG},  the  spectral data relates to the moduli space of admissible parabolic rank 2 bundles on $S/\sigma $ as seen in \cite[Section 8.3]{thesis}. 

\begin{exercise}[(*)]
Nonabelianization of the fixed point set of $\Theta_{Sp(2p,2p)}$ can also be seen through Cameral covers, as shown in \cite{ana1}. Realise the action of $\sigma$ in terms of Cameral covers. 
\end{exercise}

\begin{small}

\end{small}
\end{document}